\documentclass{amsart}
\usepackage{amssymb, amsmath}
\usepackage{mathrsfs}
\usepackage{amscd}
\usepackage{verbatim}

\allowdisplaybreaks

\theoremstyle{plain}
\newtheorem{theorem}{Theorem}[section]
\theoremstyle{remark}
\newtheorem{remark}[theorem]{Remark}
\newtheorem{example}[theorem]{Example}
\theoremstyle{plain}

\newtheorem{lemma}[theorem]{Lemma}
\newtheorem{proposition}[theorem]{Proposition}
\newtheorem{definition}[theorem]{Definition}

\numberwithin{equation}{section}

\def\R{{\mathbb R}}
\def\C{{\mathbb C}}

\newcommand{\E}{{\mathbb E}}
\renewcommand{\P}{{\mathbb P}}
\newcommand{\F}{{\mathcal F}}

\renewcommand{\a}{\alpha}

\newcommand{\g}{\gamma}

\newcommand{\e}{\varepsilon}

\renewcommand{\O}{\Omega}

\newcommand{\calL}{{\mathcal B}}
\newcommand{\n}{\Vert}

\newcommand{\embed}{\hookrightarrow}
\newcommand{\s}{^*}
\newcommand{\lb}{\langle}
\newcommand{\rb}{\rangle}

\newcommand{\limn}{\lim_{n\to\infty}}

\newcommand{\sumn}{\sum_{n\ge 1}}

\newcommand{\gL}{L^2_{\g}}

\newcommand{\VVo}{V^{0}_{\a,p}([0,T_0]\times\O;X_{a})}

\newcommand{\A}{{\mathcal A}}
\newcommand{\B}{{\mathcal B}}

\begin{document}

\title
[Damped plate and wave equation]{Structurally damped plate and
wave equations with random point force in arbitrary space dimensions}

\author{Roland Schnaubelt}
\address{Institut f\"ur Analysis \\
Universit\"at Karlsruhe (TH)\\
D-76128  Karls\-ruhe\\Germany}

\email{schnaubelt@math.uni-karlsruhe.de}

\author{Mark Veraar}
\address{Delft Institute of Applied Mathematics\\
Delft University of Technology \\ P.O. Box 5031\\ 2600 GA Delft\\The
Netherlands} \email{mark@profsonline.nl, M.C.Veraar@tudelft.nl}

\thanks{The second author was supported by the Alexander von Humboldt foundation and by a VENI subsidy 639.031.930 of the Netherlands Organization for Scientific Research (NWO)}






\date\today

\begin{abstract}
In this paper we consider structurally damped plate and wave
equations with point and distributed random forces. In order to
treat space dimensions more than one, we work in the setting of
$L^q$--spaces with (possibly small) $q\in(1,2)$. We establish
existence, uniqueness and regularity of mild and weak solutions to
the stochastic equations employing recent theory for stochastic
evolution equations in UMD Banach spaces.
\end{abstract}

\maketitle

\section{Introduction}

Structurally damped plate and wave equations have been studied intensively in
the deterministic case. In such equations the damping term has `half of the
order' of the  leading elastic term, as it has been proposed in the seminal
paper \cite{Ru86}. Point controls and feedbacks in elastic systems lead
naturally to perturbations of damped equations by Dirac measures, cf.\
\cite{LaTr}. In this paper we investigate the situation when such point
perturbations act randomly. For a one dimensional spatial domain $S$ these
problems have been treated in \cite{Mas} by means of the well-established
Hilbert space approach to stochastic partial differential equations, see e.g.\
\cite{DPZ}. However, it seems that in higher space dimensions $d\ge2$ one cannot
proceed in this way since the irregularity coming from the point measure and
the stochastic terms cannot be balanced by the smoothing effect of the analytic
semigroup for the damped plate equation. The point evaluation acts via duality
on the state space so that it becomes even more singular if one works in the
setting of $L^q(S)$ with $q>2$. Thus it is natural to look for solutions in
$L^q(S)$ with $q\in(1,2)$; in fact, we need $q\in(1,d/(d-1))$ for the plate
equation and $q\in(1,2d/(2d-1))$ for the wave equation.

In the non-damped situation the stochastic wave equation and plate equation are
studied by many authors. Here the equation is usually analyzed using Hilbert space
methods. We refer the reader to \cite{Dal,Ond,Pes} and \cite{Kimplate},
respectively, and to the references therein.

Several authors have investigated stochastic partial differential equations on
$L^q$ spaces with $q\in[2,\infty)$  (see \cite{Brz2,Kry} and the references
therein), whereas in our case $q\in(1,2)$ it seems that the only known method
is contained in the recent paper \cite{NVW3}. Our analysis is based on the
theory developed in \cite{NVW3}. Stochastic damped wave equations have been
treated in various papers during the last years, see e.g.\ \cite{BDP02a, BDPT, BrzMasSei,
CeFr, CDF, Fan, TZ}. However, it seems that random forces acting at a single
point in $S$ have been studied only in  \cite{Mas} so far.

To be concise, we will focus on the model
\begin{equation}\label{eq:plateintro}
\left\{ \begin{aligned}
&\ddot{u}(t,s) +  \Delta^2 u(t,s)-\rho \Delta \dot{u}(t,s)= f(t,s, u(t,s),\dot{u}(t,s))\\
& \qquad\qquad + b(t,s,u(t,s),\dot{u}(t,s)) \, \frac{\partial w_1(t,s)}{\partial t} +
\Big[G(t,u(t,\cdot), \dot{u}(t,\cdot))  \\
& \qquad\qquad + C(t,u(t,\cdot), \dot{u}(t,\cdot))
\frac{\partial w_2(t)}{\partial t} \Big] \delta(s-s_0),\qquad t\in [0,T],
    \ s\in S,\\
& u(0,s) = u_0(s), \ \dot{u}(0,s) = u_1(s), \qquad s\in S, \\
& u(t,s) = \Delta u(t,s) = 0,\qquad t\in[0,T], \  s\in \partial S,
\end{aligned}\right.
\end{equation}
on  a bounded domain $S\subset \R^d$ of class $C^4$. Here $\rho>0$ is a
constant and $\delta(\cdot-s_0)$ is the point mass at $s_0\in S$. The functions
$f,b,G,C$ are measurable, adapted and Lipschitz in a sense specified in
Section~\ref{sec:plate}. The process $w_1$ is a Gaussian process which is white
in time and appropriately colored in space, as discussed in
Section~\ref{sec:plate}. The process $w_2$ is a standard one--dimensional
Brownian motion which is independent of $w_1$. Note that $w_2$ drives the point
loading whereas $w_1$ governs a distributed stochastic term.

In Theorem~\ref{thm:platen2} we obtain a mild and weak solution
$(u(t),\dot{u}(t)) \in (W^{2,q}(S)\cap W^{1,q}_0(S))\times L^q(S)$ of
\eqref{eq:plateintro}, where $u$ and $\dot{u}$ possess some additional
regularity in time and in space if the initial data are regular enough. We also
state a related result for the wave equation in Theorem~\ref{thm:wave}. Our
results can be generalized in various directions. For instance, in
\eqref{eq:plateintro} one could replace the Dirichlet Laplacian by a more
general elliptic operator. One can also allow for more general nonlinearities,
see Remarks~\ref{rem:CG} and \ref{rem:nabla}, and one could treat locally
Lipschitz coefficients to some extent, see Remark~\ref{rem:local}. But for
conciseness we will focus on the setting indicated above.

For illustration, we  further give a concrete simple example for \eqref{eq:plateintro}.
\begin{example}
Consider for instance
\begin{equation}\label{eq:plateintro2}
\left\{ \begin{aligned}
&\ddot{u}(t,s) +  \Delta^2 u(t,s)-\rho \Delta \dot{u}(t,s) \\ & \qquad \qquad  = C(u(t,\cdot))
\frac{\partial w_2(t)}{\partial t} \delta(s-s_0), \ \qquad t\in [0,T],
    \ s\in S,
\\ & u(0,s) = u_0(s), \ \dot{u}(0,s) = u_1(s), \qquad s\in S, \\
& u(t,s) = \Delta u(t,s) = 0, \qquad t\in[0,T], \  s\in \partial S.
\end{aligned}\right.
\end{equation}
Here the function $C: L^1(S)\to \R$ is given by $C(x) = \int\varphi(x(s)) \, ds$, where $\varphi$ is fixed Lipschitz function $\varphi:\R\to \R$. This is a special case of \eqref{eq:plateintro} where $f,b,G$ are zero. Such a $C$ is indeed Lipschitz continuous. Moreover, if $d\leq 3$, then by Remark \ref{rem:CG} one can even take $C:C(\overline{S})\to \R$ as $C(x)= \varphi(x(s_1))$, where $s_1\in \overline{S}$ is fixed. Using Theorem \ref{thm:platen2} if $u_0$ and $u_1$ are appropriate one can obtain a unique solution $u$ to the problem \eqref{eq:plateintro2}. Moreover, $u$ satisfies a certain space-time regularity. We refer to Theorem \ref{thm:platen2} for details.
\end{example}

In Sections 2-4 we provide the necessary prerequisites for our main
results. First, we briefly discuss the theory of stochastic
integration developed in \cite{NVW1}. This theory is closely tied to
the concept of Gauss functions and operators which is also presented
in Section~\ref{sec:prel}. Based on this material, in
Section~\ref{sec:abstracteq} we recall a theorem on existence,
uniqueness and regularity of mild solutions of parabolic stochastic
equations from \cite{NVW3}. In this theorem it is possible to
consider deterministic and stochastic terms taking values in
so--called extrapolation spaces which are larger than the state
space. This fact is crucial for our approach since the Dirac
functional $\delta(\cdot-s_0)$ lives in such extrapolation spaces.
Moreover, one can use this flexibility to extend the class of
admissible processes $w_1$, see the examples in
Section~\ref{sec:plate}.

The underlying deterministic equation is studied in
Section~\ref{sec:systems}, where we  consider  the problem
\begin{equation}\label{eq:abstractwave}
\begin{split}
&\ddot{u}(t)+\rho\A^\frac{1}{2} \dot{u}(t)+\A u(t)=0, \qquad t\ge0,\\
&u(0)=u_0,\qquad  \dot{u}(0)=u_1,
\end{split}
\end{equation}
for a sectorial operator $\A$ on a Banach space $E$, see \eqref{ass:damped}.
(In \eqref{eq:plateintro}
$\A$ is the square of the Dirichlet Laplacian on $E=L^q(S)$.)
Using the operator matrix
\[
A=\begin{pmatrix} 0 & I\\ -\A & -\rho \A^\frac12
\end{pmatrix}
\qquad\text{with}\quad D(A)=D(\A)\times D(\A^\frac12)
\]
one can reformulate \eqref{eq:abstractwave} as an abstract Cauchy problem on
the state space $X=D(\A^\frac12)\times E$. It is well known that $A$ generates
an analytic semigroup on $X$ if $E=L^2(S)$, see \cite{CT89}, \cite{LaTr}.
Recently, it has been shown in \cite{CCD} that $A$ also generates an analytic
semigroup in the Banach space case. (See  \cite{CS}, \cite{FO} and \cite{JTW}
for related results.) In view of the stochastic problem,
 it is crucial to determine the  inter- and extrapolation spaces
for this semigroup, see Proposition~\ref{prop:damped}. It would be very
interesting to extend these results to damping terms which are more general
than $\rho\A^\frac{1}{2}$. In the Hilbert space case this is possible to some
extent, see \cite{CT89} and \cite{LaTr}, but this approach makes heavy use of
the Hilbert space structure. Let us explain the problem with an example, cf.\
\cite{Ru86}. In equation \eqref{eq:plateintro} it would be interesting to study
also the clamped plate equation, where $u = \frac{\partial u}{\partial n} =0$
on $\partial S$ and $n$ denotes the outer normal. If we let $\mathcal{A} =
\Delta^2$ with the above boundary conditions, then $\mathcal{A}^{\frac12}$ is
not a differential operator anymore. Instead of $\mathcal{A}^{\frac12}\dot{u}$,
we would still like to have $\Delta \dot{u}$ as a damping term, but this does
not lead to the algebraic structure of \eqref{eq:abstractwave}. Therefore, we
do not know whether the corresponding operator matrix generates a strongly
continuous semigroup if $q\neq2$.

We will write $a\lesssim b$ if there exists a universal constant $C>0$ such
that $a\leq Cb$, and $a\eqsim b$ if $a\lesssim b\lesssim a$. If the constant
$C$ is allowed to depend on some parameter $\theta$,
we write $a\lesssim_\theta b$ and $a\eqsim_\theta b$ instead.
Moreover, $X$ always denotes a Banach space,
$\B(X,Y)$ is the space of bounded linear operators from $X$ to another
Banach space $Y$, and we designate the norm in $X$ and the operator norm
by  $\|\cdot\|$.

\medskip

{\em Acknowledgment} -- The authors thank the anonymous referee for
helpful comments.

\section{Preliminaries\label{sec:prel}}
Throughout this paper  $({\O},{\F},{\P})$ is a probability space with a
filtration $(\F_t)_{t\geq 0}$. This space is used for the stochastic
equations and stochastic integrals below. In Subsection~\ref{sec:stochint}
we recall the necessary definitions and facts from the theory of stochastic
integration developed in \cite{NVW1}. As a preparation, we discuss $\gamma$-radonifying
operators and Gauss functions in the next subsection
referring to \cite{Bog,DJT,KaWe} for proofs and more details. In the last
subsection we describe a concept of Lipschitz continuity which is crucial for
our work.

\subsection{$\gamma$-radonifying operators}
In this paper,  $(\g_n)_{n\ge 1}$ always denotes a {\em Gaussian
sequence}, i.e., a sequence of independent, standard,
real-valued Gaussian random variables defined on a
 probability space $(\tilde{\O},\tilde{\F},\tilde{\P})$.
A linear operator $R:\mathcal{H}\to X$ from a separable real Hilbert space
$\mathcal{H}$ into a Banach space $X$ is called {\em a $\gamma$-radonifying operator} if for
some (and then for every) orthonormal basis $(h_n)_{n\ge 1}$ of $\mathcal{H}$
the Gaussian sum $ \sumn \g_n Rh_n$ converges in $L^2(\tilde{\O};X)$. The space
$\g(\mathcal{H},X)$ of all $\gamma$-radonifying
operators from $\mathcal{H}$ to $X$ is a
Banach space with respect to the norm
\[ \n R\n_{\g(\mathcal{H},X)} := \Big( \E \Big\n \sumn \g_n Rh_n
\Big\n^2\Big)^\frac12.\]
This norm is independent of the orthonormal basis $(h_n)_{n\ge 1}$
and the Gaussian sequence  $(\g_n)_{n\ge 1}$. It holds that
$\|R\|\le  \n R\n_{\g(\mathcal{H},X)}$. Moreover,
$\g(\mathcal{H},X)$ is an operator ideal in the
sense that if $S_1:\tilde{\mathcal{H}}\to \mathcal{H}$ and $ S_2:X\to
\tilde{X}$ are bounded operators, then $R\in \g(\mathcal{H},X)$ implies
$S_2RS_1\in \g(\tilde{\mathcal{H}},\tilde{X})$ and
\begin{equation}\label{eq:ideal}
 \n S_2RS_1\n_{ \g(\tilde{\mathcal{H}},\tilde{X})} \leq \n S_2\n \n R\n_{\g(\mathcal{H},X)} \n S_1\n.
\end{equation}
If $X$ is a Hilbert space, then $\g(\mathcal{H},X)$ is isometrically
isomorphic to  the space $\mathcal{S}^2(\mathcal{H},X)$ of Hilbert--Schmidt
operators from $\mathcal{H}$ into $X$.

We are mainly interested in the case that $\mathcal{H} = L^2(M;H)$, where $H$
is another separable real Hilbert space with inner product $[\cdot, \cdot]_H$
and $(M,\Sigma,\mu)$ is a $\sigma$-finite measure space. Let $\Phi:M\to
\B(H,X)$. Assume that $\Phi\s x\s \in L^2(M;H)$ for all $x\s\in X\s$ and that
there exists an $R\in \g(L^2(M;H),X)$ such that
\[ \lb Rf,x\s\rb = \int_M [f(t), \Phi\s(t)x\s]_H\,d\mu(t)\]
for all $f\in L^2(M;H)$ and $x\s\in X\s$.
Then we say that $R$ is {\em represented} by $\Phi$. In this case
$\Phi$ is called a \emph{Gauss function}, and we write
$\Phi\in\g(M;H,E)$ and
\[\|\Phi\|_{\g(M;H,X)}:= \|R\|_{\g(L^2(M;H),X))}.\]
We write $\g(M;X)$ instead of $\g(M;\R,X)$. If there is no danger
of confusion we will identify $R$ and $\Phi$, cf.\ Subsection~2.3 in
\cite{NVW1}. For a Hilbert space $X$, we have $\g(M;H,X) =
L^2(M;\mathcal{S}^2(H,X))$ isometrically.

For the space $X=L^p(S)$, $p\in [1, \infty)$, the following
well-known square function estimate gives a useful way to verify that
$\Phi:M\to \B(H,X)$ is a Gauss function, see \cite[Proposition~2.1]{BP99} and \cite[Theorem 2.3]{BrzvN03} (also see \cite[Proposition~6.1]{NVWco}):
\begin{equation}\label{eq:square}
\|\Phi\|_{\g(M;H,X)} \eqsim_p \Big\|\Big(\int_M \sum_{n\geq 1} |\Phi(t) h_n|^2
\, d\mu(t)\Big)^{\frac12}\Big\|_{L^p(S)}.
\end{equation}

\subsection{Stochastic integration in UMD spaces\label{sec:stochint}}
We now discuss the stochastic integral for processes $\Phi:[0,T]\times\O\to
\calL(H,X)$ as it  was introduced and investigated in \cite{NVW1}. Here $X$ is
a UMD Banach space and $H$ is a separable real Hilbert space. The reader is
referred to \cite{Bu3} and \cite{NVW1} concerning UMD spaces. But, for the
present paper, it suffices to recall that the reflexive $L^q$, Sobelev,
Bessel--potential and Besov spaces are UMD spaces.

We denote by $L^0(\Omega;E)$ the vector space of all equivalence classes
of measurable functions from $\Omega$ to a Banach space $E$, and we endow
$L^0(\Omega;E)$ with the convergence in probability.
A process $\Phi:[0,T]\times\O\to \calL(H,X)$ is called {\em $H$--strongly
measurable} if $\Phi h$ is strongly measurable in $X$ for all $h\in H$, where
we let $(\Phi h)(t,\omega):= \Phi(t,\omega)h$. The process $\Phi$ is called
\emph{$H$--strongly adapted} if the map $\omega\mapsto \Phi(t,\omega) h$ is
$\F_t$--strongly measurable  for all $t\in [0,T]$ and $h\in H$.
We also set $\Phi_\omega(t)=\Phi(t,\omega)$.

An {\em $H$-cylindrical Brownian motion} is a family
$W_H=(W_H(t))_{t\in [0,T]}$ of bounded linear operators from $H$ to
$L^2(\O)$ satisfying
\begin{enumerate}
\item $W_H h = (W_H(t)h)_{t\in [0,T]}$ is a real-valued Brownian motion for
          each $h\in H$,
\item $ \E (W_H(s)g \cdot W_H(t)h) = (s\wedge t)\,[g,h]_{H}$ for all
$s,t\in [0,T]$ and  $g,h\in H.$
\end{enumerate}

Let  $0\le a<b<T$, $A\subset \O$ be $\F_a$-measurable, $x\in X,$
and  $h\in H$.
The stochastic integral of the indicator process $1_{(a,b]\times A}\otimes
(h\otimes x)$ is then defined as
\[\int_0^T 1_{(a,b]\times A}\otimes (h\otimes x)\,dW_H
     := 1_A(W_H(b)h - W_H(a)h) x.\]
(Analogously, one defines the integral  for the trivial process $1_{[0]\times
A}\otimes (h\otimes x)$.) By linearity, this definition extends to adapted step
processes $\Phi:[0,T]\times\O\to \calL(H,X)$ whose values are finite rank
operators.
 An $H$-strongly measurable and adapted process $\Phi$ is called {\em
stochastically integrable} with respect to $W_H$ if there exists a sequence of
adapted step processes $\Phi_n:[0,T]\times \O \to\calL(H,X)$ with values in the
finite rank operators from $H$ to $X$ and a pathwise continuous process
$\zeta:[0,T]\times\O\to X$ such that the following two conditions are
satisfied:
\begin{enumerate}
\item $\limn \langle \Phi_n h, x^*\rangle  = \langle\Phi h, x^*\rangle$
in measure on $[0,T]\times\Omega$ for all $h\in H$ and $x^*\in X^*$;
\item $\displaystyle \limn \int_0^\cdot \Phi_n(t)\,dW_H(t) = \zeta$ \quad in
$L^0(\O;C([0,T];X))$.
\end{enumerate}
In this situation,
$\zeta$ is uniquely determined as an element of $L^0(\O;C([0,T];X))$ and
it is
called the {\em stochastic integral} of $\Phi$ with respect to $W_H$. We write
\[ \zeta = \int_0^\cdot \Phi\,dW_H = \int_0^\cdot \Phi(t)\,dW_H(t).\]
The process $\zeta$ is a continuous local martingale starting at
zero, see \cite[Theorem~5.5]{NVW1}.

\begin{proposition}\cite[Theorems 5.9 and 5.12]{NVW1}\label{prop:NVW}
Let $X$ be a UMD space. For an $H$-strongly measurable and adapted process
$\Phi:[0,T]\times\O\to\calL(H,X)$ the following assertions are equivalent.
\begin{enumerate}
\item The process $\Phi$ is stochastically integrable with respect to $W_H$.

\item For all $x\s\in X\s$ the process $\Phi\s x\s$ belongs to
$L^0(\O;L^2(0,T;H))$, and there exists a pathwise continuous process $\zeta\in
L^0 (\O;C([0,T], X))$ such that for all $x\s\in X\s$ we have
\[\lb \zeta, x\s\rb =  \int_0^\cdot \Phi\s x\s\,dW_H \qquad \hbox{ in
$L^0(\O;C([0,T]))$;}\]

\item $\Phi_\omega\in \g(0,T;H,X)$ for a.e.\ $\omega\in\Omega$.
\end{enumerate}

In this situation we have $\zeta = \int_0^{\cdot} \Phi\,dW_H$ in
$L^0(\O;C([0,T];X))$. Furthermore, for all $p\in (1, \infty)$,
\[
\E \sup_{t\in [0,T]}\Big\n \int_0^t \Phi\,dW_H\Big\n^p \eqsim_{p,X}
\E\n \Phi\n_{\g(0,T;H,X)}^p.
\]
\end{proposition}

\subsection{$L^2_{\g}$-Lipschitz functions}
We now treat a class of Lipschitz functions which is needed in the existence
result for the stochastic equation presented in the next section. See
\cite{NVW3} for more details.

Let $(M,\Sigma)$ be a countably generated measurable space and let $\mu$ be a
finite measure on $(M,\mu)$. Then $L^2(M,\mu)$ is separable. We then define
 \[L^2_{\g}(M,\mu;X) :=
\g(M,\mu;X)\cap L^2(M,\mu;X),\] which is a Banach space endowed with the norm
\[\|\phi\|_{L^2_{\g}(M,\mu;X)} :=
  \|\phi\|_{\g(M,\mu;X)} + \|\phi\|_{L^2(M,\mu;X)}.\]
Note that the simple functions are dense in $L^2_{\g}(M,\mu;X)$.

Let $H$ be a separable real Hilbert space, let $X_1$ and $X_2$ be Banach
spaces, and let $f:M\times X_1\to \calL(H,X_2)$ be a function such that for all
$x\in X_1$ we have $f(\cdot, x)\in \g(L^2(M,\mu;H),X_2)$. For simple functions
$\phi:M\to X_1$ one easily checks that the map $s\mapsto f(s, \phi(s))$
belongs to $\g(L^2(M,\mu;H),X_2)$. We call $f$ an {\em $\gL$--Lipschitz function with
respect to $\mu$} if $f$ is strongly continuous in the second variable and
we have
\begin{equation}\label{eq:glipsimp}
\|f(\cdot,\phi_1) - f(\cdot, \phi_2)\|_{\g(L^2(M,\mu;H),X_2)}\le
C\|\phi_1-\phi_2\|_{L^2_{\g}(M,\mu;X_1)}
\end{equation}
for a constant $C\ge0$ and
all simple functions $\phi_1, \phi_2:M\to X_1$.
In this case the mapping $\phi\mapsto f(\cdot,\phi(\cdot))$ extends uniquely to
a Lipschitz mapping from $L^2_{\g}(M,\mu;X_1)$ into $\g(L^2(M,\mu;H),X_2)$. Its
Lipschitz constant will be denoted by $L^{\g}_{\mu,f}$. Finally,
if $f$ is $\gL$-Lipschitz with respect to all finite measures $\mu$ on
$(M,\Sigma)$ and
\[L^\g_{f} := \sup\{L^{\g}_{\mu,f}: \mu \ \text{is a finite measure on
$(M,\Sigma)$}\}\] is finite, then we say that $f$ is a {\em
$\gL$--Lipschitz function}.

In the next lemma we state a
simpler sufficient condition for the  $\gL$--Lipschitz property.
However, for this result one has to impose
an additional restriction on the Banach space $X_2$ which we first
 introduce. Let  $p\in [1,2]$, and let $(r_j)_{j\geq 1}$ be a Rademacher
sequence,  i.e., $(r_{j})_{j\geq 1}$ is an
independent, identically  distributed sequence with
$\P(r_1 = 1) = \P(r_1 = -1)=\frac12$. A Banach
space $X$ has {\em type $p$} if there exists a constant
$C_p\ge 0$ such that for all $x_1,\dots,x_n\in E$ we have
\begin{equation}\label{type}
\Big(\E \Big\n \sum_{j=1}^n r_j\, x_j\Big\n^2\Big)^\frac12\le C_p
\Big(\sum_{j=1}^n \n x_j\n^p\Big)^\frac1p.
\end{equation}
For more information on this concept we refer the reader to \cite{DJT}
and the references
therein. We recall that every Banach space has type $1$, that the
spaces $L^p(S)$, $1\le p<\infty$, have type $\min\{p,2\}$  and that
Hilbert spaces have type $2$. The property has a certain ordering:
If $X$ has type $p$, then $X$ has type $\tilde{p}$ for all $1\leq
\tilde{p}<p$ as well. Furthermore, every UMD space has nontrivial
type, i.e., type $p$ for some $p\in (1,2]$. But we will not need
this fact. In type 2 spaces the $\gL$--Lipschitz property can be
checked using only the norm in $\g(H,X_2)$.

\begin{lemma}\cite[Lemma~5.2]{NVW3}\label{lem:type2Lipschitz}
Let $X_2$ be a space with type $2$. Let $f:M\times X_1\to \g(H,X_2)$ be a
function such that $f(\cdot,x)$ is strongly measurable for each $x\in X_1$. If
there is a constant $C$ such that
\begin{align}
\label{eq:type2bdd} \|f(s,x)\|_{\g(H,X_2)} & \le  C(1+\|x\|),\\
\label{eq:type2Lipschitz} \|f(s,x) - f(s,y)\|_{\g(H,X_2)}&\le C\|x-y\|
\end{align}
for all $s\in M$ and $x,y\in X_1$,
then $f$ is a $\gL$--Lipschitz function.
We also have $L_{f}^\g\le C_2 C$ for the constant $C_2$ from \eqref{type}.
Moreover, $f$ satisfies
\[\|f(\cdot,\phi)\|_{\g(L^2(M,\mu;H),X_2)} \le C_2 C(1+\|\phi\|_{L^2(M,\mu;X_1)}).\]
\end{lemma}
If $f$ does not depend on $M$, one can check that
\eqref{eq:glipsimp} implies \eqref{eq:type2bdd} and
\eqref{eq:type2Lipschitz}.
Clearly, every $\gL$-Lipschitz function $f:X_1\to \g(H,X_2)$ is a Lipschitz
function. The converse does not hold (see \cite[Theorem~1]{NVgammaLip}).
The next example shows that standard
substitution operators are $\gL$--Lipschitz.

\begin{example}\cite[Example~5.5]{NVW3}
Let $p\in [1, \infty)$,  $(M,\Sigma, \mu)$ be a finite measure space, and
$b:\R\to \R$ be  Lipschitz continuous. Define the Nemytskii map $B:L^p(M)\to
L^p(M)$ by $B(\varphi)(s) := b(\varphi(s)).$ Then $B$ is
$\gL$--Lipschitz with respect to $\mu$.
\end{example}

\section{The abstract stochastic evolution equation\label{sec:abstracteq}}
Recall that  $(\O,\F,\P)$ is a probability space with filtration
 $(\F_t)_{t\in[0,T]}$. Let $H_1$ and $H_2$ be separable real Hilbert spaces.
Let $X$ be a UMD Banach space and let $Y$ be a Banach space.
On the Banach space $X$ we consider the problem
\begin{equation}\tag{SE}\label{eq:SEtype}
\left\{\begin{aligned}
dU(t) & = (A U(t) + F(t,U(t)) + \Lambda_G G(t,U(t)))\,dt  \\ & \qquad + B(t,U(t))\,dW_{H_1}(t) + \Lambda_C C(t,U(t)) \, d W_{H_2}(t), \  t\in [0,T],\\
 U(0) & = U_0.
\end{aligned}
\right.
\end{equation}
We assume that $A$ generates an analytic $C_0$--semigroup
$(S(t))_{t\geq 0}$ on $X$. Thus, there are constants $M\ge1$ and $w_0\in\R$
such that $\|S(t)\|\leq  M e^{w_0 t}$ for $t\geq 0$.  Let $w>w_0$.
Our further
assumptions make use of the fractional power scale associated to $A$,
see e.g.\ \cite{Am}. For $a \in [0,1]$,
we define the space $X_a = D((w-A)^a)$
with the norm $\|x\|_a=\|(w-A)^a x\|$.  For  $\theta \in [0,1]$,
we further introduce the extrapolation space $X_{-\theta}$
which is the completion of $X$ with respect to the norm
$\|x\|_{-\theta} = \|(w-A)^{-\theta} x\|$.
 The operator $A$  has a restriction (extension) to an operator
on the space $X_a$ (the space $X_{-\theta}$) which generates the analytic
$C_0$--semigroup given by the restrictions (extensions) of $S(t)$ on the space
$X_a$ (the space $X_{-\theta}$). We usually denote the restrictions and
extensions again by $A$ and $S(t)$. Moreover, $(w-A)^\beta$ is an isomorphism
from $X_\alpha$ to $X_{\alpha-\beta}$, where $-1\le \alpha-\beta \le \alpha\le
1$. Finally, $X_\alpha$ is continuously embedded into $X_{\alpha-\beta}$.

Going back to (SE), we now list the assumptions on the linear operators
$\Lambda_j:Y\to X_{-\theta_j}$ for $j=G,C$ and  on the functions
\begin{align*}
&F:[0,T]\times \Omega\times X_a\to X,
 & &G:[0,T]\times \Omega\times X_a\to Y,\\
&B:[0,T]\times\O\times X_a\to \calL(H_1,X_{-\theta_B}),
 & &C:[0,T]\times \O\times X_a\to \calL(H_2,Y).
\end{align*}
Here the exponents $a,\theta_G,\theta_B,\theta_C$ belong to $[0,1]$, but
in the next theorem we impose further restrictions.
Moreover, the initial value $U_0:\Omega\to X_a$ has to be
strongly $\F_0$--measurable. The interval $[0,T]$ is endowed with
the Borel $\sigma$--algebra $\mathcal{B}_{[0,T]}$.

\begin{enumerate}
\let\ALTERWERTA\theenumi
\let\ALTERWERTB\labelenumi
\def\theenumi{(H1)}
\def\labelenumi{(H1)}
\item \label{as:analytic}
$A$ generates an analytic strongly continuous
semigroup $(S(t))_{t\geq 0}$ on $X$.
\let\theenumi\ALTERWERTA
\let\labelenumi\ALTERWERTB
\let\ALTERWERTA\theenumi
\let\ALTERWERTB\labelenumi
\def\theenumi{(H2)}
\def\labelenumi{(H2)}
\item \label{as:LipschitzFtype}
The map $(t, \omega)\mapsto F(t, \omega,x)\in X$ is strongly measurable and adapted
for each $x\in X_a$. The function
$F$ has linear growth and is Lipschitz continuous in $x$ uniformly in
$[0,T]\times\O$; i.e., there are constants $L_F,C_F\ge0$ such that
\begin{align*}
\|F(t,\omega,x)-F(t,\omega,y)\|_{X}&\leq  L_F\|x-y\|_a,
\\ \nonumber
\|F(t,\omega,x)\|_{X}&\leq  C_F(1+\|x\|_a)
\end{align*}
for all $t\in [0,T]$, $\omega\in \O$, and $x,y\in X_{a}$.
\let\theenumi\ALTERWERTA
\let\labelenumi\ALTERWERTB
\let\ALTERWERTA\theenumi
\let\ALTERWERTB\labelenumi
\def\theenumi{(H3)}
\def\labelenumi{(H3)}
\item \label{as:LipschitzGtype}
The map
$(t, \omega)\mapsto G(t, \omega,x)\in Y$ is strongly measurable and adapted
for all $x\in X_a$. The function $\Lambda_G G$ has linear growth and
is Lipschitz continuous in $x$
uniformly in $[0,T]\times\O$; i.e., there are constants $L_G,C_G\ge0$ such
that
\begin{align*}
\|\Lambda_G (G(t,\omega,x)-G(t,\omega,y))\|_{-\theta_G}&\leq
L_G\|x-y\|_a,
\\ \nonumber
\|\Lambda_G G(t,\omega,x)\|_{-\theta_G}&\leq  C_F(1+\|x\|_a)
\end{align*}
for all $t\in [0,T]$, $\omega\in \O$, and $x,y\in X_{a}$.
\let\theenumi\ALTERWERTA
\let\labelenumi\ALTERWERTB
\let\ALTERWERTA\theenumi
\let\ALTERWERTB\labelenumi
\def\theenumi{(H4)}
\def\labelenumi{(H4)}
\item \label{as:LipschitzBtype}
The map
$(t,\omega)\mapsto B(t,\omega, x)\in \calL(H_1,X_{-\theta_B})$
is $H_1$-strongly measurable and adapted  for all $x\in X_a$. The function $B$
is $\gL$-Lipschitz of
linear growth uniformly in $\O$; i.e., there are constants $L_B^{\g}$ and
$C_B^{\g}$ such that
\begin{align*}
\qquad\quad  \| B(\cdot,\omega,\phi_1) -
B(\cdot,\omega,\phi_2)\|_{\g((0,T),\mu;H_1,X_{-\theta_B})}
  &\le L_B^{\g}\|\phi_1-\phi_2\|_{L^2_{\g}((0,T),\mu;X_a)},\\
\|B(\cdot,\omega,\phi_1)\|_{\g((0,T),\mu;H_1,X_{-\theta_B})}
&\le C_B^{\g}(1+\|\phi_1\|_{L^2_{\g}((0,T),\mu;X_a)}).
\end{align*}
 for all finite measures $\mu$ on $([0,T],
\mathcal{B}_{[0,T]})$,  for all $\omega\in \O$, and  all $\phi_1, \phi_2\in
L^2_{\g}((0,T),\mu;X_a)$.
\let\theenumi\ALTERWERTA
\let\labelenumi\ALTERWERTB
\let\ALTERWERTA\theenumi
\let\ALTERWERTB\labelenumi
\def\theenumi{(H5)}
\def\labelenumi{(H5)}
\item \label{as:LipschitzCtype}
The map
$(t,\omega)\mapsto \Lambda_C C(t,\omega, x)\in
\calL(H_2,X_{-\theta_C})$ is $H_2$-strongly measurable and
adapted  for all $x\in X_a$. The composition $\Lambda_C C$ is
$\gL$-Lipschitz of linear growth uniformly in $\O$; i.e., there are constants
$L_C^{\g}$ and $C_C^{\g}$ such that
\begin{align*}
\qquad\;\; \| \Lambda_C ( C(\cdot,\omega,\phi_1) -
C(\cdot,\omega,\phi_2))\|_{\g((0,T),\mu;H_2,X_{-\theta_C})}
& \le L_C^{\g}\|\phi_1-\phi_2\|_{L^2_{\g}((0,T),\mu;X_a)}\\
\|\Lambda_C C(\cdot,\omega,\phi_1)\|_{\g((0,T),\mu;H_2,X_{-\theta_C})}
  &\le C_C^{\g}(1+\|\phi_1\|_{L^2_{\g}((0,T),\mu;X_a)})
\end{align*}
for all finite measures $\mu$ on $([0,T],
\mathcal{B}_{[0,T]})$, for  all $\omega\in \O$, and all $\phi_1, \phi_2\in
L^2_{\g}((0,T),\mu;X_a)$.
\let\theenumi\ALTERWERTA
\let\labelenumi\ALTERWERTB
\end{enumerate}

For $p\in [1, \infty)$ and  $\alpha\in (0,\frac12)$ we define
$V^{0}_{\a,p}([0,T]\times\O;X)$ as the linear space of continuous
adapted processes $\phi:[0,T]\times \O\to X$ such that
\[\|\phi(\cdot,\omega)\|_{C([0,T];X)}+ \Big(\int_0^T \|s\mapsto (t-s)^{-\alpha} \phi(s,\omega)\|_{\g(L^2(0,t),X)}^p \,dt\Big)^{\frac1p} <\infty\]
for almost all $\omega\in\Omega$.
In $V^{0}_{\a,p}([0,T]\times\O;X)$ we identify indistinguishable processes;
i.e., processes $\phi_1$ and $\phi_2$ such that a.s.\ for all $t\in [0,T]$ we have
$\phi_1(t) = \phi_2(t)$.

In order to introduce our solution concept, we
recall some notation from \cite{NVW3}.
For $\phi\in L^1(0,T;X_{-\theta})$ with $\theta\in [0,1)$, we  write
\begin{equation}\label{eq:detconv}
S*\phi(t) = \int_0^t S(t-s) \phi(s) \, ds, \qquad t\in [0,T].
\end{equation}
Young's inequality and the regularity properties of $S(t)$ yield
$S*\phi\in L^1(0,T;X)$.
For $j=1, 2$ and processes $\Phi:[0,T]\times\O\to \calL(H_j,X_{-\theta})$ with $\theta\in
[0,\frac12)$ which are $H_j$-strongly measurable and adapted and such that for all $t\in
[0,T]$  the map
\[s\mapsto S(t-s)\Phi(s) \quad \text{belongs to} \quad \g(0,t;H_j,X),\]
almost surely, we set
\begin{equation}\label{eq:stochconv1}
S\diamond_j \Phi(t) = \int_0^t S(t-s) \Phi(s) \, d W_{H_j}(s).
\end{equation}
This integral exists for each $t\in [0,T]$ due to Proposition~\ref{prop:NVW}.

\begin{definition}\label{def:mild}
An $X_a$--valued process $(U(t))_{t\in [0,T]}$ is called a {\em mild solution} of
\eqref{eq:SEtype} if
\renewcommand{\labelenumi}{(\roman{enumi})}
\renewcommand{\theenumi}{(\roman{enumi})}
\begin{enumerate}
\item $U:[0,T]\times\O\to X_a$ is strongly measurable and adapted,

\item $F(\cdot, U)\in L^0(\O;L^1(0,T;X))$,

\item $\theta_G\in [0, 1)$ and $\Lambda_G G(\cdot, U) \in
L^0(\O;L^1(0,T;X_{-\theta_G}))$,

\item\label{it:mildB} for all $t\in [0,T]$, $(s,\omega)\mapsto S(t-s)B(s,U(s))$
is $H_1$-strongly measurable and adapted and belongs to $\g(0,t;H_1,X)$ almost
surely,

\item\label{it:mildC} for all $t\in [0,T]$, $(s,\omega)\mapsto S(t-s)\Lambda_C
C(s,U(s))$ is $H_2$-strongly measurable and adapted and belongs to
$\g(0,t;H_2,X)$ almost surely,

\item for all $t\in [0,T]$, the following equality holds a.s.\ in $X$:
\[U(t) = S(t) U_0 + S*F(\cdot,U)(t)+S*\Lambda_G G(\cdot,U)(t) + S\diamond_1 B(\cdot,U)(t) +  S\diamond_2 \Lambda_C C(\cdot,U)(t).\]
\end{enumerate}
\renewcommand{\labelenumi}{(\arabic{enumi})}
\renewcommand{\theenumi}{(\arabic{enumi})}
\end{definition}

\begin{definition}\label{def:weak}
An $X_a$--valued process $(U(t))_{t\in [0,T]}$ is called a {\em weak
solution} of (SE) if
\renewcommand{\labelenumi}{(\roman{enumi})}
\renewcommand{\theenumi}{(\roman{enumi})}
\begin{enumerate}
\item\label{en:weak2} $U$ is strongly measurable and adapted and has paths in $L^1(0,T;X_a))$ a.s.,

\item\label{it:weakF} $F(\cdot, U)\in L^0(\O;L^1(0,T;X))$,

\item\label{it:weakG} $\theta_G\in [0, 1)$ and $\Lambda_G G(\cdot, U)
\in L^0(\O;L^1(0,T;X_{-\theta_G}))$,

\item\label{it:weakB} $\theta_B\in [0,\frac12)$ and $B(\cdot, U):
[0,T]\times\O\to \calL(H_1,X_{-\theta_B})$ is $H_1$-strongly measurable with
\[\int_0^T \|B(t,U(t))\|_{\calL(H_1,X_{-\theta_B})}^2\, dt <\infty \ \text{almost surely,}\]

\item\label{it:weakC} $\theta_C\in [0,\frac12)$ and $\Lambda_C C(\cdot, U)$ is
$H_2$-strongly measurable with
\[\int_0^T \|\Lambda_C C(t,U(t))\|_{\calL(H_2,X_{-\theta_C})}^2\, dt <\infty \ \text{almost surely,}\]

\item\label{en:weak4} for all $t\in [0,T]$ and $x^*\in D(A^*)$, we have
\begin{equation}\label{weaksol}\begin{aligned}
\lb  U(t), &x^*\rb - \lb u_0, x^*\rb   \\
& = \int_0^t \lb U(s),A^*x^*\rb
 + \lb F(s,U(s)) + \Lambda_G G(s,U(s)), x^*\rb \, ds \\
&\qquad  + \int_0^t B^*(s,U(s)) x^* \, d W_{H_1}(s)
         + \int_0^t (\Lambda_C C(s,U(s)))^* x^* \, d W_{H_2}(s),
\end{aligned}
\end{equation}
almost surely.
\end{enumerate}
\renewcommand{\labelenumi}{(\arabic{enumi})}
\renewcommand{\theenumi}{(\arabic{enumi})}
\end{definition}

The following result and its proof are standard in stochastic
evolution equations (cf. \cite[Theorem 5.4]{DPZ}).
Since our setting slightly differs from the existing literature, we
include a short proof.

\begin{proposition}\label{prop:weakmild}
Let $X$ be a UMD space. Let $\theta_G\in [0,1)$,
$\theta_B,\theta_C\in [0,\frac12)$ and let $a\in [0,\frac12)$. Then the
following assertions hold.
\begin{enumerate}
\item If $U\in L^0(\O;L^1(0,T;X_a))$ is a mild solution of (SE) such that
Definition \ref{def:weak} \ref{it:weakF}--\ref{it:weakC} hold, then
$U$ is a weak solution of (SE).

\item If $U$ is a weak solution of (SE) such that Definition \ref{def:mild}
\ref{it:mildB}--\ref{it:mildC} hold, then $U$ is a mild solution of (SE).
\end{enumerate}
\end{proposition}

\begin{proof}
We will write $\tilde{F} = F+\Lambda_G G$, $H = H_1\times H_2$, and
$\tilde{B} = (B,\Lambda_C C)$.

(1) \ Let $t\in [0,T]$ and $x^*\in D(A^*)$. From the
definition of a mild solution, Proposition~\ref{prop:NVW} and the
(stochastic) Fubini theorem we obtain that almost surely
\begin{align*}
&\int_0^t  \lb U(s), A^*  x^*\rb\,ds \\
& = \int_0^t \lb  u_0, S(s)^*A^*x^*\rb\,ds
  + \int_0^t \int_r^t \lb \tilde{F}(r,U(r)), S(s-r)^*A^*x^*\rb \,ds\,dr\\
& \qquad +\int_0^t \int_r^t \tilde{B}^*(r,U(r)) S^*(s-r) A^*x^* \,ds\,dW_H(r)\\
& = \lb S(t)u_0,x^*\rb -\lb u_0, x^* \rb
    +\int_0^t \lb S(t-r) \tilde{F}(r,U(r)), x^* \rb \, dr
    - \int_0^t \lb \tilde{F}(r, U(r)), x^* \rb \, dr\\
&\qquad  + \int_0^t \tilde{B}^*(r,U(r)) S^*(t-r) x^* \, d W_H(r)
      - \int_0^t \tilde{B}^*(r,U(r)) x^*\, d W_H(r)\\
& = \lb U(t), x^*\rb - \lb u_0, x^*\rb
      - \int_0^t  \lb \tilde{F}(r,U(r)), x^*\rb \, dr
- \int_0^t \tilde{B}^*(r,U(r)) x^* \, d W_H(r).
\end{align*}
This shows that $U$ is a weak solution.

(2) \ Fix $t\in [0,T]$. Let $f\in C^1([0,t])$, $x^*\in D(A^*)$, $\varphi
= f\otimes x^*$, and $U$ be a mild solution. It\^o's formula implies that
\begin{align}
\label{eq:varsol} \lb U(t), \varphi(t)\rb  & = \lb u_0,
\varphi(0)\rb  + \int_0^t \lb U(s), A^* \varphi(s)\rb  + \lb
\tilde{F}(s,U(s)), \varphi(s)\rb \, ds
\\ \nonumber & \qquad + \int_0^t \lb U(s), \varphi'(s)\rb \, ds + \int_0^t \tilde{B}^*(s,U(s)) \varphi(s) \, d
W_H(s),
\end{align}
almost surely. By linearity one can extend \eqref{eq:varsol} to functions
$\varphi:[0,t]\to D(A^*)$ of the form $\varphi = \sum_{n=1}^N
f_n\otimes x^*_n$, with $f_n\in C^1([0,t])$ and $x_n^*\in D(A^*)$
for all $n=1, \ldots, N$. By density this extends to all $\varphi\in
C^1([0,t];D(A^*))$. In particular, for $x^*\in D((A^*)^2)$ we can take
$\varphi(s) = S^*(t-s) x^*$ and thus deduce
\[\begin{aligned}
  \lb U(t), x^*\rb - \lb S(t) u_0, x^*\rb
& =  \int_0^t \lb S(t-s) \tilde{F}(s, U(s)), x^*\rb \, ds\\
& \qquad + \int_0^t \tilde{B}^*(s, U(s)) S^*(t-s) x^* \, d W_H(s),
\end{aligned}\]
almost surely. Since the integrals in Definition~\ref{def:mild} exist
by our assumptions, the Hahn-Banach theorem yields
that $U$ is a mild solution.
\end{proof}

We can now formulate the main abstract existence and uniqueness result which is
a consequence of  Theorems~7.1 and 7.2 in \cite{NVW3}.

\begin{theorem}\label{thm:mainexistenceLloc}
Let $X$ be a UMD space with type $\tau\in [1,
2]$ and suppose that \ref{as:analytic}-\ref{as:LipschitzCtype} are
satisfied. Assume that $0\le a+\theta_G<\frac32 -\frac1\tau$ and
$a+\max\{\theta_B, \theta_C\}<\frac12$. Let $U_0:\O\to X_a$ be
strongly $\F_0$--measurable. Then the following assertions hold.
\begin{enumerate}
\item\label{part:existunique} If $\alpha\in (0,\frac12)$ and $p>2$ are such
that $a+\max\{\theta_B,\theta_C\}<\alpha-\frac1p$, then there exists a unique
mild solution $U\in \VVo$ of \eqref{eq:SEtype}.

\item\label{part:regularity} Let $\lambda\ge 0$ and $\delta\geq a$ satisfy
$\lambda+\delta<\min\{1-\theta_G, \frac12-\theta_B,
\frac12-\theta_C\}$. Then the mild solution $U$ of \eqref{eq:SEtype}
has a version such that almost all paths satisfy $U-SU_0\in
C^\lambda([0,T];X_{\delta})$.
\end{enumerate}
\end{theorem}
\begin{proof}
In \cite{NVW3} the problem
\begin{equation}\tag{SE$'$}\label{eq:SEtype2}
\left\{\begin{aligned}
dU(t) & = (A U(t) + \tilde{F}(t,U(t))) \,dt + \tilde{B}(t,U(t))\,dW_{H}(t), \  t\in [0,T],\\
 U(0) & = U_0
\end{aligned}
\right.
\end{equation}
has been considered. Clearly, \eqref{eq:SEtype} can be written as
\eqref{eq:SEtype2} if we take $\tilde{F} = F+\Lambda_G G$, $H = H_1\times H_2$
and $\tilde{B} = (B,\Lambda_C C)$. In this way the result follows
immediately from Theorems~7.1 and 7.2 in \cite{NVW3}.
\end{proof}

\begin{remark}
There is a version of Theorem~\ref{thm:mainexistenceLloc} for
locally Lipschitz functions as well (see \cite[Section 8]{NVW3}).
Some of the results below remain true for locally Lipschitz
coefficients. However, for the sake of simplicity  we concentrate on
the (global) Lipschitz case here.
\end{remark}

\section{Strongly damped second order equations\label{sec:systems}}
Before we turn to the equation \eqref{eq:plateintro}, we have to treat
a class of deterministic damped second order equations. We investigate
the problem
\begin{equation}\label{eq:damped}
\begin{split}
&\ddot{u}(t)+\rho\A^\frac{1}{2} \dot{u}(t)+\A u(t)=0, \qquad t\ge0,\\
&u(0)=u_0,\qquad  \dot{u}(0)=u_1,
\end{split}
\end{equation}
and, for $\alpha\in (\frac12,1]$, its variant
\begin{equation}\label{eq:damped2}
\begin{split}
&\ddot{u}(t) +\A^\alpha(\rho \dot{u}(t)+\A^{1-\alpha} u(t))=0, \qquad t\ge0,\\
&u(0)=u_0,\qquad  \dot{u}(0)=u_1,
\end{split}
\end{equation}
both on a Banach space $E$ with norm $\|\cdot\|_0$. We assume that
\begin{equation}\label{ass:damped}
\begin{split}
&\A \text{ is invertible on $E$, } \overline{D(\A)}=E, \ \lambda\in \rho(-\A)
\text{ and } \|\lambda (\lambda I+\A)^{-1}\|_{\B(E)}\le M\\
&\text{for all } \lambda\in\C\setminus \{0\} \text{ with }
|\arg\lambda|\le \pi-\phi
\text{ and some } \phi\in\Big(0,\frac{\pi}{2}\Big), \; M>0.\\
&\text{Further, let \ }\alpha\in\Big[\frac{1}{2},1\Big], \;\;\rho>0, \quad
\text{ and } \rho>2\cos \frac{\pi-\phi}{2} \text{ \ \ if } \alpha=\frac{1}{2}.
\end{split}\end{equation}
We denote by $E_\theta$ with $\theta\in [-1,1]$ the domains of fractional powers for
$\A$ on $E$; i.e., for $\theta\in [0,1]$ we take $E_{\theta} = D(\A^{\theta})$, and
for $\theta\in [-1,0)$ we let $E_{\theta}$ be the completion of $E$ with respect
to the norm $\|x\|_{\theta} = \|\A^{\theta}x\|$.
Concerning \eqref{eq:damped} we look for solutions
$u\in C^2(\R_+,E)\cap C^1(\R_+,E_\frac{1}{2}) \cap C(\R_+,E_1)$, whereas
the solutions of \eqref{eq:damped2} have to satisfy
$u\in C^2(\R_+,E)\cap C^1(\R_+,E_\frac{1}{2})$
and $\rho \dot{u}+\A^{1-\alpha} u\in C(\R_+,E_\alpha)$.

In the recent paper \cite{CCD}, it was shown that
the operator matrix
\begin{equation}\label{eq:A-damped}
\begin{split}
A&=\begin{pmatrix}
0 & I\\ -\A & -\rho \A^\alpha
\end{pmatrix} \quad \text{with the domain} \\
D(A)&= \{(\varphi,\psi)\in E_{\frac{3}{2}-\alpha}\times E_{\frac{1}{2}}:
\A^{1-\alpha}\varphi+\rho\psi\in E_\alpha\}
\end{split}\end{equation}
generates an analytic $C_0$--semigroup on $X=E_{\frac{1}{2}}\times E$, where the
action of $A$ is defined by $A(\varphi,\psi)= (\psi,
-\A^\alpha(\A^{1-\alpha}\varphi+2\rho\psi))$ if $\alpha>1/2.$ In the most
important case $\alpha=1/2$, we simply obtain $D(A)=E_1\times E_{\frac{1}{2}}$
and the matrix in \eqref{eq:A-damped} is understood in the usual way. In
\eqref{ass:damped} the constant $\rho>0$  has to satisfy an additional lower
bound if $\alpha=\frac12$. This restriction cannot be avoided in view of
Remark~1.1 of \cite{CCD}. However, in the typical applications (as those
discussed below) one can choose $\phi>0$ arbitrarily small, so that
\eqref{ass:damped} holds for all $\rho>0$ in these applications. We recall that
due to a result of H\"ormander \cite{Ho}, $A$ does not generate a
$C_0$-semigroup if, say, $E=L^p(\R^d)$, $p\neq2$, $\rho=0$, and $\A$ is the
negative Laplacian.
 Moreover, in the Hilbert space case and for a strictly
positive self adjoint $\A$, Chen and Triggiani proved that $A$ generates an
analytic semigroup if one replaces the damping term $\rho \A^\alpha$ by a self
adjoint operator $B$ satisfying $\rho_1 \A^\alpha\le B\le \rho_2 \A^\alpha$
in form sense for some $\rho_2>\rho_1>0$ and $\alpha\in[1/2,1]$, see
\cite{CT89} and \cite{LaTr}. It would be very interesting to extend this result
to the Banach
space setting. In \cite{CT89} it was also shown that for $\alpha<1/2$, $A$ does
not generate an analytic semigroup.

In the next result we use the real interpolation spaces
for a sectorial operator $C$ on Banach space $Y$, see e.g.\
\cite{Lun}, \cite{Tr1}. Recall that
\begin{equation}\label{eq:int-emb}
\begin{aligned}
(Y,D(C))_{\gamma+\e,q}  \embed (Y,D(C))_{\gamma,1} &\embed D((w-C)^\gamma)
\embed \ldots
\\ & \ldots \embed (Y,D(C))_{\gamma,\infty} \embed (Y,D(C))_{\gamma-\e,q}
\end{aligned}
\end{equation}
for every
$q\in[1,\infty]$ and $0<\gamma-\e\le\gamma<\gamma+\e<1$. Moreover, if
$C_{-1}:X\to X_{-1}$ is the extrapolation of $C$, then $X_{\gamma-1}$ is
isomorphic to the domain  $D((w - C_{-1})^\gamma)$ in $X_{-1}$, see
\cite[Theorem~V.1.3.8]{Am}. So the isomorphism \eqref{eq:inter/extra-damped2}
below implies that
\begin{equation}\label{eq:inter/extra-damped2b}
E_{\frac{1}{2}-(1-\alpha)(\theta-\e)} \times E_{-\alpha(\theta-\e)}
\embed  X_{-\theta} \embed
E_{\frac{1}{2}-(1-\alpha)(\theta+\e)} \times E_{-\alpha(\theta+\e)},
\end{equation}
for all $\e>0$ with $0<\theta-\e<\theta<\theta+\e\le 1/2$, and analogous embeddings hold in
case of \eqref{eq:inter/extra-damped3} and \eqref{eq:inter/extra-damped4b}.
We write $X\cong Y$ if $X$ and $Y$ are canonically isomorphic.

\begin{proposition}\label{prop:damped} Assume that
\eqref{ass:damped} holds on a Banach space $E$. Then the operator
matrix $A$ from \eqref{eq:A-damped} generates an analytic $C_0$--semigroup on
$X=E_{\frac{1}{2}}\times E$. For all $(u_0,u_1)\in D(A)$ the
problems \eqref{eq:damped} and \eqref{eq:damped2} have unique
solutions in the above specified sense. Moreover,  we have
\begin{align}\label{eq:inter/extra-damped1a}
X_\theta&= E_{\frac{1}{2}+(1-\alpha)\theta}\times E_{\alpha\theta},
\\
 X_{-\frac{1}{2}} & \cong E_\frac{\alpha}{2}\times  E_{-\frac{\alpha}{2}}\,,
 \label{eq:inter/extra-damped1b} \\
(X_{-1},X)_{1-\theta,q} &\cong (E,E_1)_{\frac{1}{2}-(1-\alpha)\theta,q}
                           \times (E_{-1},E)_{1-\alpha\theta,q}
\label{eq:inter/extra-damped2}
\end{align}
for all $\theta\in[0,1/2]$ and $q\in[1,\infty]$.

If, additionally $\alpha=1/2$, then
\begin{align}\label{eq:inter/extra-damped3}
(X, D(A))_{\theta,q} &= (E,E_1)_{\frac{1+\theta}{2},q}
                           \times (E,E_1)_{\frac{\theta}{2},q}\,,\\
X_{-1}&\cong E\times  E_{-\frac{1}{2}}\,, \label{eq:inter/extra-damped4a}
\\ (X_{-1},X)_{1-\theta,q} &\cong
(E,E_1)_{\frac{1-\theta}{2},q}
          \times (E_{-1},E)_{1-\frac{\theta}{2},q}
\label{eq:inter/extra-damped4b}
\end{align}
for all $\theta\in(0,1)$ and $q\in[1,\infty]$. Furthermore, if $E$ is
reflexive, then  $X^* = (E^*)_{-\frac12} \times E^*$ and
\begin{equation}\label{eq:adjoint}
A^*=\begin{pmatrix} 0 & -\A^*
\\ I & -\rho (\A^*)^\frac12
\end{pmatrix}
\qquad\text{with}\quad D(A^*) = E^* \times (E^*)_{\frac12}\,,
\end{equation}
where $(E^*)_\theta$ is the fractional power space for $\A^*$.
\end{proposition}

\begin{proof}
The generation property was shown in \cite[Theorem~2.3]{CCD}. It easily
implies the unique solvability of \eqref{eq:damped} and \eqref{eq:damped2}. The
equation \eqref{eq:inter/extra-damped1a} was also  proved in
\cite[Theorem~2.3]{CCD}. (We note that in \cite{CCD} it was assumed that $E$
is reflexive.  However, this property is not needed in the parts of the
proofs which are relevant to us.)
Take $(\varphi,\psi)\in X_\frac{1}{2}=
 E_{1-\frac{\alpha}{2}}\times E_{\frac{\alpha}{2}}$. Due to
\cite[p.2316]{CCD}, we have
\[A^{-1}= \begin{pmatrix} -\rho\A^{\alpha-1} &-\A^{-1}\\I & 0\end{pmatrix}.\]
Using \eqref{eq:inter/extra-damped1a}, we can estimate
\[ \|A^{-1}(\varphi,\psi)\|_{X_\frac{1}{2}}
  \eqsim \|\A^{1-\frac{\alpha}{2}}(\rho\A^{\alpha-1}\varphi + \A^{-1}\psi)\|_0
         +\|\varphi\|_\frac{\alpha}{2}
\lesssim_{\rho} \|\varphi\|_{\frac{\alpha}{2}}+\|\psi\|_{-\frac{\alpha}{2}}\,,
\]
where $\|\cdot\|_0$ denotes the norm on $E$. Conversely, we obtain
\begin{align*}
\|\varphi\|_{\frac{\alpha}{2}} + \|\psi\|_{-\frac{\alpha}{2}}
&=\|\varphi\|_{\frac{\alpha}{2}} + \|\A^{-\frac{\alpha}{2}}\psi
 +\rho \A^{\frac{\alpha}{2}}\varphi - \rho \A^{\frac{\alpha}{2}}\varphi\|_0\\
&\lesssim_{\rho} \|\varphi\|_{\frac{\alpha}{2}}
 +\|\A^{1-\frac{\alpha}{2}}(\rho\A^{\alpha-1}\varphi + \A^{-1}\psi)\|_0
\eqsim \|A^{-1}(\varphi,\psi)\|_{X_{\frac{1}{2}}}.
\end{align*}
The isomorphism \eqref{eq:inter/extra-damped1b} thus follows since
$X_{-\frac{1}{2}}$ is isomorphic to the completion of $X_{\frac{1}{2}}$ with
respect to the norm $\|A^{-1}(\varphi,\psi)\|_{X_\frac{1}{2}}$, cf.\
\cite[Theorem~V.1.3.8]{Am}. Notice that real interpolation respects cartesian
products due to its definition via the $K$--functional. Furthermore,
the reiteration theorem (see e.g.\ \cite[Theorem~1.2.15]{Lun}) implies that
$(X_{-1},X)_{1-\theta, q} = (X_{-\frac12},X)_{1-2\theta, q}$. The
equality \eqref{eq:inter/extra-damped2} is then a consequence
of \eqref{eq:inter/extra-damped1b} and reiteration.

Let $\alpha=1/2$. In this case we have $X_1=E_1 \times E_\frac{1}{2}$. Take
$(\phi,\psi)\in X$. We first show \eqref{eq:inter/extra-damped3}. We estimate
as above
\begin{align*}
\|A^{-1}(\varphi,\psi)\|_X &
  \eqsim \|\A^\frac{1}{2}(\rho\A^{-\frac{1}{2}}\varphi+ \A^{-1}\psi)\|_0
         +\|\varphi\|_0
     \lesssim_{\rho}  \|\varphi\|_0 + \|\psi\|_{-\frac{1}{2}},\\
\|\varphi\|_0 + \|\psi\|_{-\frac{1}{2}}
&= \|\A^\frac{1}{2}( \A^{-1}\psi+ \rho\A^{-\frac{1}{2}}\varphi
      -\rho\A^{-\frac{1}{2}}\varphi )\|_0  +\|\varphi\|_0
\lesssim_{\rho}  \|A^{-1}(\varphi,\psi)\|_X.
\end{align*}
The formulas \eqref{eq:inter/extra-damped3} and \eqref{eq:inter/extra-damped4b}
can now be established as the isomorphism \eqref{eq:inter/extra-damped2}.

The last assertion follows easily from $(E_{\frac12})^* =
(E^*)_{-\frac12}$ (see Theorem~V.1.4.12 of \cite{Am}) and a
straightforward calculation using that the operator matrix in
\eqref{eq:adjoint} is invertible in $X^*$.
\end{proof}

\section{The stochastically perturbed damped plate equation\label{sec:plate}}
In this section we prove existence, uniqueness and regularity results for the
structurally damped plate equation with noise, given by
\begin{equation}\label{eq:plate}
\left\{
\begin{aligned}
&\ddot{u}(t,s) + \Delta^2 u(t,s) - \rho \Delta \dot{u}(t,s)= f(t,s, u(t,s),\dot{u}(t,s)) \\
& \qquad \qquad + b(t,s,u(t,s),\dot{u}(t,s)) \,
     \frac{\partial w_1(t,s)}{\partial t}
         + \Big[G(t,u(t,\cdot), \dot{u}(t,\cdot)) \\
&\qquad \qquad + C(t,u(t,\cdot), \dot{u}(t,\cdot))
\frac{\partial w_2(t)}{\partial t} \Big] \delta(s-s_0), \qquad
     t\in [0,T],\ s\in S,\\
& u(0,s) = u_0(s), \ \dot{u}(0,s) = u_1(s), \qquad s\in S,
\\ & u(t,s) = \Delta u(t,s) = 0,\qquad  t\in[0,T], \ s\in \partial S,
\end{aligned}\right.
\end{equation}
Using Proposition~\ref{prop:damped} and the theory from Section~\ref{sec:abstracteq},
we will reformulate
problem \eqref{eq:plate} as an equation of the type \eqref{eq:SEtype}.
 We first list our assumptions
and notations, where  subsets $M$ of $\R^n$ are endowed with
the Borel $\sigma$--algebra $\B_M$.
\begin{itemize}
\item[(A0)] Let $S\subset \R^d$ be a bounded domain with boundary $\partial S$
of class $C^4$ and  $(\O,\F,\P)$ be a probability space with filtration
$(\F_t)_{t\geq 0}$.  The number $\rho>0$ and the point $s_0\in S$ are fixed,
and $\delta$ denotes the usual point mass.
\end{itemize}

\noindent
Let $H_1$ be a separable real Hilbert space, $H_2=Y=\R$ and
$E=L^q(S)$ for some $q\in (1,\infty)$. We identify $Y$  with
$\calL(H_2, Y)$,
where  we interpret each $y\in Y$ as the operator $h\mapsto y h$.
We further introduce the negative Dirichlet Laplacian on $E$ by
\[\mathcal{B} \varphi = -\Delta \varphi, \quad
 D(\mathcal{B})  = W^{2,q}(S)\cap W^{1,q}_0(S).\]
We set $\A = \B^2$, so that $\A^{\frac12} = \B$. As in
Section~\ref{sec:systems} we define  the operator $(A,D(A))$ on $X
:= E_{\frac12}\times E$  by setting
\[A = \left(%
\begin{array}{cc}
  0 & I \\
  -\mathcal{A} & -\rho \mathcal{A}^{\frac12} \\
\end{array}%
\right), \qquad D(A) = E_{1}\times E_{\frac12}.
\]
Since $\A$ is a sectorial operator of angle $\phi$ for all $\phi\in
(0,\pi/2)$ (see e.g.\ \cite[Theorem~8.2]{DHP}), the assumption
\eqref{ass:damped} is satisfied for the above $\A$ and $\rho>0$. So
Hypothesis~\ref{as:analytic}  in Section~3 follows from
Proposition~\ref{prop:damped}, i.e., $A$ generates an analytic
$C_0$--semigroup $(S(t))_{t\ge0}$. We also  recall that
\begin{equation}\label{eq:DBtheta}
D(\B^\theta) =\begin{cases}
   H^{2\theta,q}(S), \quad \text{if \ }  0\le 2\theta<\frac1q,\\
 \{\varphi\in H^{2\theta,q}(S): \varphi=0 \text{ \ on } \partial S \},
  \qquad \text{if \ }  \frac1q<2\theta\le 1.
\end{cases}
\end{equation}
(cf.\ \cite[Corollary 2.2]{Am00}).
We also  observe that equation \eqref{eq:inter/extra-damped1a}
with $\alpha=\frac12$ gives
\begin{equation}\label{eq:fractdomain}
X_{\delta} = E_{\frac12 +\frac12 \delta}\times E_{\frac12 \delta}
\end{equation}
for $\delta\in [0,\frac12]$.
Combining this identity with  \eqref{eq:DBtheta}, we deduce
\[  X_{\delta}     =\begin{cases}
          (H^{2+2\delta,q}(S)\cap W^{1,q}_0(S)) \times H^{2\delta,q}(S),
              \quad \text{if }2\delta\in(0, \frac{1}{q}),\\
   \{(\varphi,\psi)\in H^{2+2\delta,q}(S) \times H^{2\delta,q}(S):
        \varphi=\Delta \varphi=\psi = 0 \text{ on } \partial S\}, \;\;
         \text{if  }2\delta\in( \frac{1}{q},1).
\end{cases}
\]
We further make the following hypotheses.
\begin{itemize}
\item[(A1)] The functions $f,b:[0,T]\times\O\times S\times \R\times \R\to \R$
are jointly measurable, adapted to $(\F_t)_{t\geq 0}$, and
 Lipschitz functions and of linear growth in the fourth and fifth variable,
uniformly in the other variables. The process $w_2$ is a
standard real--valued Brownian motion with respect to $(\F_t)_{t\geq 0}$.
We set $W_{H_2}(t):=w_2(t)$ for $t\ge0$.
\item[(A2)] The maps $G,C:[0,T]\times\O\times X\to \R$ are jointly
measurable, adapted  to $(\F_t)_{t\geq 0}$, and
Lipschitz and of linear growth in the third variable,
uniformly in the other variables.
\item[(A3)] The process $w_1$ can be written in the form
$i_1 W_{H_1}$, where $i_1\in\B(H_1,L^r(S))$ for some $r\in [1,\infty]$
and $W_{H_1}$  is a cylindrical Wiener process
 with respect to $(\F_t)_{t\geq 0}$ being independent of $W_{H_2}$.
\end{itemize}

\noindent
Assumption (A3) has to be interpreted in the sense that
\[w_1(t,s) = \sum_{n\geq 1} (i_1 h_n)(s) W_{H_1}(t) h_n,
\qquad  t\in \R_+, s\in S,\] where $(h_n)_{n\geq 1}$ is an orthonormal basis
for $H_1$ and the sum converges in $L^r(S)$. In the examples below we will be
more specific about $i_1$ and $W_{H_1}$. Typically, $H_1$ is $L^2(S)$ or the
reproducing kernel Hilbert space and $i_1 i_1^*$ is the covariance operator of
$w_1$, see e.g.\ \cite{Bog}, \cite{BrzvNe} and the references therein. It is
also possible to assume that $i_1$ takes values in an extrapolation space of
$L^r(S)$, but we do not consider this generalization here.

We now continue to reformulate \eqref{eq:plate} as
a problem of the type \eqref{eq:SEtype}. We focus on the case
$a=0$ in (H2)--(H5), though  we comment on possible extensions
in some remarks below. We define $F:[0,T]\times\O\times X\to X$ by
\[F(t,\omega,x)(s) = \left(%
\begin{array}{c}
  0 \\
  f(t,\omega, s, x(s),\dot{x}(s)) \\
\end{array}%
\right).\] It straightforward to check that $F$
satisfies~\ref{as:LipschitzFtype} because of (A1). Let $\theta_B\in
[0,\frac12)$. The map   $B:[0,T]\times\O\times X\to
\calL(H_1,X_{-\theta_B})$ is defined by
\begin{equation}\label{eq:b-def}
B(t,\omega,x) h (s) = \left(%
\begin{array}{c}
  0 \\
  b(t,\omega, s, x(s),\dot{x}(s)) (i_1 h)(s) \\
\end{array}%
\right).
\end{equation}
It will be assumed that $B$ satisfies \ref{as:LipschitzBtype}.
Below, we discuss various classes of examples where this assumption
holds. We further take a suitable $\theta_G=\theta_C\in (0,\frac12)$
and define $\Lambda_C=\Lambda_G=\Lambda\in \B(\R, X _{-\theta_C})$
by
\[\Lambda y (s) = \left(%
\begin{array}{c}
  0 \\
  \delta(s-s_0) y \\
\end{array}%
\right).
\]

We claim that
for all $1<q<\frac{d}{d-1}$ if $d\geq 2$ and all $1<q<\infty$ if $d=1$, there
exists a $\theta_C\in (\frac{d}{2 q'},\frac12)$ such that
$\Lambda$ is well-defined. Indeed, let $\theta<\theta_C$.  Due to equation
\eqref{eq:inter/extra-damped2b} with $\alpha=\frac12$, it holds
\[E_{\frac12-\frac12\theta} \times E_{-\frac12\theta} \hookrightarrow
D((-A)^{-\theta_C}) = X_{-\theta_C}.
\]
So we have to find a $\theta\in [0,\frac12)$ with
 $\delta(\cdot-s_0)\in E_{-\frac12\theta}$. Theorem~V.1.4.12 of \cite{Am}
implies that $E_{-\frac12\theta}=(D((\B^*)^{\theta}))^*$. It thus remains to
show that the point evaluation $\delta_{s_0}:D((\B^*)^{\theta})\to \R$,
$\delta_{s_0} f = f(s_0)$, defines a bounded linear map. Since $\B^*$ is the
realization of the  negative Dirichlet Laplacian  on $L^{q'}(S)$, we deduce
from \eqref{eq:DBtheta}
that $D((\B^*)^{\theta})$ is a closed subset of $H^{2\theta,q'}(S)$.
Sobolev's embedding theorem (cf.\ \cite[Theorem~4.6.1]{Tr1}) yields
$H^{2\theta,q'}(S)\hookrightarrow C(\overline{S})$ if
$2\theta>\frac{d}{q'}$. So the claim follows.

Assertion (A2) now implies \ref{as:LipschitzGtype}
and \ref{as:LipschitzCtype} since  $\Lambda G$ and  $\Lambda C$
factorize through the spaces $Y=\R$ and $\B(H_2,Y)=\R$, respectively.
(Use the ideal property \eqref{eq:ideal}.)

Summing up, we have found spaces $X,Y,H_1,H_2$, maps $A,F,\Lambda_G,G,
B,\Lambda_C,C$, and processes $W_{H_1}, W_{H_2}$ for which we can formulate the
equation (SE) from Section~3. In view of Theorem~\ref{thm:mainexistenceLloc},
this problem  has a unique mild solution $U$ which we call a \emph{mild
solution} of \eqref{eq:plate}.

To justify this notion of a mild solution to \eqref{eq:plate},
we need to define a weak solution of \eqref{eq:plate}. To this aim,
we assume
that $D(\mathcal{A}^*)\hookrightarrow C(\overline{S})$. One easily
checks that this embedding always holds for $d=1, 2, 3,4$ and for all
$q<\frac{d}{d-4}$ if $d\geq 5$. Assume that $f,G,b,C$ are as before.
We say that a process $u:[0,T]\times\O\times S\to \R$ is a {\em weak solution}
of \eqref{eq:plate} if it is measurable,
$u(t,\cdot)$ is $\F_t\otimes \mathcal{B}_S$-measurable  for all $t\in [0,T]$,
$u\in W^{1,2}(0,T;L^q(S))$ a.s., and for all $\phi\in W^{2,q'}(S)\cap
W^{1,q'}_0(S) = D(\mathcal{A}^*)$ we have
\begin{align}
\nonumber &\lb  u(t,\cdot), \phi\rb -  \lb u_0, \phi\rb - t \lb u_1,\phi\rb
 +\int_0^t\int_0^{\sigma} \lb u(\tau,\cdot), \Delta^2 \phi\rb\,d\tau\,d\sigma\\
\nonumber &\qquad - \rho \int_0^t \lb u(\sigma,\cdot), \Delta \phi\rb\,d\sigma
     + t\rho \lb u_0, \Delta\phi\rb \\
\label{eq:weakplate} & = \int_0^t \int_0^{\sigma} \Big( \lb f(\tau,\cdot, u(\tau,\cdot),\dot{u}(\tau,\cdot)),
\phi\rb +  G(\tau,u(\tau,\cdot), \dot{u}(\tau,\cdot)) \phi(s_0) \Big)\,
d\tau \, d\sigma
\\ \nonumber& \qquad + \int_0^t \int_0^{\sigma} \lb b(\tau,\cdot, u(\tau,\cdot),\dot{u}(\tau,\cdot)),
\phi\rb \, dw_1(\tau) \, d\sigma
\\ \nonumber& \qquad + \int_0^t \int_0^{\sigma} C(\tau,u(\tau,\cdot), \dot{u}(\tau,\cdot)) \phi(s_0) \,
dw_2(\tau) \, d\sigma,
\end{align}
where $\lb \cdot, \cdot\rb$ denotes the $(L^q(S),L^{q'}(S))$-duality.

We will show below that \eqref{eq:plate} has a unique mild $U$ and a
unique weak solution $u$ satisfying $U = (u, \dot{u})$. In the next
theorem we use the notation introduced above.

\begin{theorem}\label{thm:platen2}
Let $d\geq 1$ and $1<q<\frac{d}{d-1}$. Assume  that (A0)--(A3)
hold,  that $B$ satisfies \ref{as:LipschitzBtype}, and that $u_0:\O\to
W^{2,q}(S)\cap W^{1,q}_0(S)$ and $u_1:\O\to L^q(S)$ are
$\F_0$-measurable. Then the following assertions hold.
\begin{enumerate}
\item For all $\alpha\in (0,\frac12)$ and $p>2$ with
$\max\{\theta_B,\theta_C\}<\alpha-\frac1p$, there exists a unique mild solution
$U$ of \eqref{eq:plate} belonging to
\begin{equation}\label{eq:solspace}
V^{0}_{\a,p}([0,T]\times\O;(W^{2,q}(S)\cap W^{1,q}_0(S)) \times L^q(S)).
\end{equation}
There is a unique weak solution $u\in W^{1,2}(0,T;L^q(S))$ of
\eqref{eq:plate} such that $(u,\dot{u})$ belongs to the space in
\eqref{eq:solspace}. Moreover, $U = (u, \dot{u})$.

\item There exists a version of $u$ with paths that satisfy $u\in
C([0,T];W^{2,q}(S)\cap W^{1,q}_0(S))$ and $\dot{u} \in C([0,T];L^q(S))$. \item
Let $\eta\in (0,\frac12]$. If $u_0:\Omega\to E_{\frac12 +\frac12\eta}$ and
$u_1:\Omega\to E_{\frac12\eta}$, then there exists a version of $u$ with paths
that satisfy $u\in C^{\lambda}([0,T];E_{\frac12 +\frac12 \delta})$ and $\dot{u}
\in C^{\lambda}([0,T];E_{\frac12 \delta})$ for all $\delta,\lambda\geq 0$ with
$\delta+\lambda<\min\{\eta, \frac12-\theta_B, \frac12-\theta_C\}$.
\end{enumerate}
\end{theorem}

Note that if $d=1$, then one can take $q\in (1, \infty)$ arbitrary, and when restricting to $q\in [2, \infty)$ one could also apply the theory from \cite{Brz2} to obtain the above result. However, if $d\geq 2$, we need that $q<2$ and therefore require the theory from \cite{NVW3}.
\begin{proof}
We have already formulated \eqref{eq:plate} as \eqref{eq:SEtype}. Let
$\alpha\in (0,\frac12)$ and $p>2$ be such that
$\max\{\theta_B,\theta_C\}<\alpha-\frac1p$. Set $U_0=(u_0,u_1)$.
Theorem~\ref{thm:mainexistenceLloc}\eqref{part:existunique} gives a unique mild
solution $U\in V^{0}_{\a,p}([0,T]\times\O;X)$ of \eqref{eq:SEtype2}, where we
set $\tilde{F} = F+\Lambda G$, $H = H_1\times H_2$ and $\tilde{B} = (B,\Lambda
C)$.  It is given by
\begin{equation}\label{mild-sol}
U(t)=S(t)U_0 + \int_0^t S(t-s)\tilde{F}(s,U(s))\,ds + \int_0^t S(t-s)\tilde{B}(s,U(s))\,dW_H(s)
\end{equation}
almost surely. Write $U = (u, v)$. We show that $\dot{u} = v$ and
that $u$ is a weak solution of \eqref{eq:plate}. Fix $t\in [0,T]$.
Tt follows from Proposition~\ref{prop:weakmild} that $U$ is a weak
solution of \eqref{eq:SEtype2}. Hence, for all $x^*\in D(A^*)$, we have
\begin{equation}\label{eq:weaktilde}
\begin{aligned}
\lb U(t), x^*\rb - \lb U_0, x^*\rb  = \int_0^t \lb U(s), A^*x^*\rb & + \lb
\tilde{F}(s,U(s)), x^*\rb \, ds \\ & + \int_0^t \tilde{B}^*(s,U(s)) x^* \, d
W_{H}(s),
\end{aligned}
\end{equation}
almost surely. In particular, for $x^* = (\phi, 0)$ with $\phi\in
E^*$, the equations \eqref{eq:weaktilde} and \eqref{eq:adjoint} yield
that $\lb u(t,\cdot), \phi\rb - \lb u_0, \phi\rb
= \int_0^t \lb v(\tau,\cdot), \phi \rb \, d\tau$ almost
surely. Therefore, $\dot{u} = v$ almost surely. Moreover, if we take
$x^* = (0,\phi)$ with $\phi\in D(A^*)$ and use \eqref{eq:adjoint}
again, we obtain
\begin{align}
\label{eq:weakalmost} \lb & \dot{u}(t,\cdot), \phi\rb -  \lb u_1, \phi\rb +
\int_0^{t} \lb u(\tau,\cdot), \Delta^2 \phi\rb \, d\tau
- \rho\int_0^t\lb \dot{u}(\tau,\cdot), \Delta \phi\rb\,d\tau
\\ \nonumber & = \int_0^t \Big(\lb f(\tau,\cdot, u(\tau,\cdot),\dot{u}(\tau,\cdot)),
\phi\rb +  G(\tau,u(\tau,\cdot), \dot{u}(\tau,\cdot)) \phi(s_0) \Big)\, d\tau
\\ \nonumber  & \quad + \int_0^t \lb b(\tau,\cdot, u(\tau,\cdot),\dot{u}(\tau,\cdot)),
\phi\rb \, dw_1(\tau)
\\ \nonumber & \quad + \int_0^t C(\tau,u(\tau,\cdot),
\dot{u}(\tau,\cdot)) \phi(s_0) \, dw_2(\tau)
\end{align}
almost surely. Now integration with respect to $t$  yields the result.

To show that $u$ is the unique weak solution, we show that every
weak solution gives a mild solution $U=(u, \dot{u})$.
The assumptions yield $u(0, \cdot) = u_0$  and $\dot{u}(0, \cdot) = u_1$
in $L^q(S)$. Fix $t\in [0,T]$ and $\phi\in D(\A^*)$.
Equation \eqref{eq:weakalmost} follows from \eqref{eq:weakplate} by
differentiation  with respect to $t$.
We claim that \eqref{eq:weaktilde} holds for all $x^*\in D(A^*)$.
For $x^*= (\phi, 0)$ with $\phi\in E^*$ this is clear from
\eqref{eq:adjoint} and $u(t,\cdot) - u_0 = \int_0^t
\dot{u}(\tau,\cdot) \, d\tau$ almost surely. For $x^* = (0,\phi)$
with $\phi\in D(\mathcal{A}^*)$ one can check that \eqref{eq:weaktilde}
reduces to \eqref{eq:weakalmost}, using \eqref{eq:adjoint} again.
By linearity and density we obtain
\eqref{eq:weaktilde} for all $x^*\in D(A^*)$. Now Proposition
\ref{prop:weakmild} implies that $U$ is a mild solution of
\eqref{eq:plate}.

Theorem~\ref{thm:mainexistenceLloc}\eqref{part:regularity}  shows that
$U-S (u_0,u_1)$ has paths in $C^{\lambda}([0,T];D((-A)^{\delta}))$
for all $\delta,\lambda\geq 0$ with
$\lambda+\delta<\min\{1-\theta_G, \frac12-\theta_B, \frac12-\theta_C\}$.
By the assumption in (3) and equation
\eqref{eq:fractdomain} we have $(u_0,u_1)\in D((-A)^\eta)$.
Therefore, $S (u_0,u_1)\in C^{\lambda}([0,T];D((-A)^{\delta}))$ a.s.\ whenever
$\lambda+\delta<\eta$. Now the result follows from \eqref{eq:fractdomain}.
\end{proof}

\begin{remark}\label{rem:CG}
We indicate an extension of the above result to the case where
$C,G:[0,T]\times\O\times C(\overline{S})\times E\to \R$ if $d\le 3$.
(Observe that in this case one can allow for point evaluations
in the third coordinate of $C$ and $G$.) First, we
 note that the identity \eqref{eq:fractdomain} yields
\[X_a\hookrightarrow (W^{2+2\tilde{a},q}(S) \times W^{2\tilde{a},q}(S))\]
for all $\tilde{a}\in[0,a)$, where we must have  $a\in [0,\frac12)$ in view of
Theorem~\ref{thm:mainexistenceLloc}. Sobolev's embedding leads to
$W^{2+2\tilde{a},q}(S)\hookrightarrow C(\overline{S})$  if
\begin{equation}\label{atilde}
2+2\tilde{a}-\frac{d}{q}>0 \iff q>\frac{d}{2+2\tilde{a}}\,.
\end{equation}
If \eqref{atilde} holds for some
$0\le \tilde{a}<a<\frac{1}{2}$ and $q\in(1,d/(d-1))$, then there is a version
of Theorem~\ref{thm:platen2} which is valid for $C$ and $G$ defined
only for $\phi\in C(\overline{S})$ (provided that $a<\min\{\theta_B,\theta_C, \theta_G\}-\frac12$).
For $d=1$ and $d=2$ the condition \eqref{atilde}
holds even for $a=\tilde{a}=0$ and all $q>1$.
For $d=3$ and each $1<q <3/2= d/(d-1)$,
 we can find an arbitrarily small $\tilde{a}>0$ satisfying \eqref{atilde}.
Therefore we can choose $a\in (\tilde{a},\frac12-\theta_B)$ if (H4) holds
for some $\theta_B<\frac12$.
For $d\geq 4$, the inequality \eqref{atilde} contradicts $q<d/(d-1)$
and $\tilde{a}<1/2$.
\end{remark}

We now discuss the interplay between $b$ and $w_1$ in several examples, where
we specify $H_1$, $i_1$ and $\theta_B$. Throughout the examples below $W_{H_1}$
is a cylindrical Brownian process as in (A3). We start with the case when the
Brownian motion is colored in space.
\begin{example}\label{ex:infty}
Assume that the covariance $Q_1\in \calL(L^2(S))$ of $w_1$ is compact. Then
there exist numbers $(\lambda_n)_{n\geq 1}$ in $\R_+$ and an orthonormal system
$(e_n)_{n\geq 1}$ in $L^2(S)$ such that
\begin{equation}\label{eq:Q1}
Q_1= \sum_{n\geq 1} \lambda_n e_n \otimes e_n
\end{equation}
Assume that
\begin{equation}\label{eq:assQ-inf}
\sum_{n\geq 1} \lambda_n \|e_n\|_{\infty}^2<\infty.
\end{equation}
Let $H_1=L^2(S)$ and $i_1:L^2(S)\to L^\infty(S)$ be given by $i_1 = \sum_{n\geq
1} \sqrt{\lambda_n} e_n \otimes e_n$. Then \ref{as:LipschitzBtype} is satisfied
with $a=\theta_B = 0$.
\end{example}

It will be clear from the proof that \eqref{eq:assQ-inf} can be replaced by
\begin{equation}\label{eq:assQ}
\Big(\sum_{n\geq 1} \lambda_n |e_n|^2\Big)^{\frac12}\in L^\infty(S).
\end{equation}

\begin{remark}
A symmetric and positive operator $Q\in\calL(L^2(S))$ maps $L^2(S)$
continuously into $L^\infty(S)$ if only if \eqref{eq:Q1} and \eqref{eq:assQ}
hold. Indeed, if $Q$ satisfies \eqref{eq:Q1} and \eqref{eq:assQ}, then the
Cauchy-Schwarz inequality implies that
\[|\sqrt{Q} h(s)| = \Big|\sum_{n\geq 1} \sqrt{\lambda_n} e_n(s) [e_n,h]_{L^2(S)} \Big|\leq \Big(\sum_{n\geq 1} \lambda_n |e_n(s)|^2\Big)^{\frac12} \|h\|_{L^2(S)}\]
for almost all $s\in S$ and all $h\in  L^2(S).$
Conversely, if $\sqrt{Q}:L^2(S)\to L^\infty(S)$ is bounded, then it is well-known that $\sqrt{Q}\in
\calL(L^2(S))$ is Hilbert-Schmidt and therefore compact.
In particular, there exists an orthonormal basis $(e_n)_{n\geq 1}$ in
$L^2(S)$ and $(\lambda_n)_{n\geq 1}$ in $\R_+$ such that \eqref{eq:Q1} holds.
Now, for almost all $s\in S$ we estimate
\[\begin{aligned}
\Big(\sum_{n\geq 1} \lambda_n |e_n(s)|^2\Big)^{\frac12} &=
\sup_{\|\beta\|_{\ell^2}\leq 1} \Big|\sum_{n\geq 1} \sqrt{\lambda_n} e_n(s)
\beta_n \Big| \\ & = \sup_{\|\beta\|_{\ell^2}\leq 1}\Big|\sqrt{Q}
\Big(\sum_{n\geq 1} \beta_n e_n\Big)(s)\Big|\leq \|\sqrt{Q}\|_{\calL(L^2(S),
L^\infty(S))}.
\end{aligned}\]
\end{remark}

\begin{proof}[Proof of Example~\ref{ex:infty}]
Our assumptions imply that  $i_1\in \calL(L^2(S),L^\infty(S))$.
Equation \eqref{eq:b-def} thus defines a map
$B:[0,T]\times\O\times X\to \calL(H_1,X)$. Moreover, the function
$(t,\omega)\mapsto B(t,\omega, x)$ is $H_1$-strongly measurable and adapted
in $X$, for each $x\in X$.  We check the $\gL$-Lipschitz property. Let $\mu$
be a finite measure on $[0,T]$. We have  to show that
\begin{align*}
\| B(\cdot,\omega,\phi_1) -  B(\cdot,\omega,\phi_2)\|_{\g(L^2((0,T),\mu;H,X)}
   &\le C\, \|\phi_1-\phi_2\|_{L^2_{\g}((0,T),\mu;X)}
\end{align*}
for all $\phi_1, \phi_2\in L^2_{\g}((0,T),\mu;X)$ and some constants $C\ge0$.
We write $\phi_1 = (\phi_{11},\phi_{12})$ and $\phi_2 =(\phi_{21},\phi_{22})$ with
 $\phi_{i1}\in \gL((0,T),\mu;E_{\frac12})$ and
 $\phi_{i2}\in \gL((0,T), \mu;E)$ for $i=1, 2$. Recall that formula
\eqref{eq:square} says that
\[\|\Phi\|_{\g(L^2((0,T),\mu;H),L^q(S))}\eqsim_q
\Big\|\Big( \int_0^T \sum_{n\geq 1} |\Phi(t) e_n|^2 \,
d\mu(t)\Big)^{\frac12}\Big\|_{L^q(S)}\] for all $\Phi\in
\g(L^2((0,T),\mu;H),L^q(S))$. Using that $b$ is Lipschitz, we thus obtain
\begin{align*}
\| B(\cdot,\omega,\phi_1)& - B(\cdot,\omega,\phi_2)\|_{\g(L^2((0,T),\mu;H),X)}\\
&\eqsim_q \Big\|\Big( \int_0^T |b(t,\omega, \phi_{11},\phi_{12})
   - b(t,\omega, \phi_{21}, \phi_{22})|^2 \, d\mu(t)
     \sum_{n\geq 1} |i_1  e_n|^2\Big)^{\frac12}\Big\|_{E}\\
 & \lesssim_{Q_1} \Big\|\Big( \int_0^T |b(t,\omega, \phi_{11},\phi_{12})
 - b(t,\omega, \phi_{21}, \phi_{22})|^2 \, d\mu(t) \Big)^{\frac12}\Big\|_{E}\\
& \lesssim L_b \Big\|\Big( \int_0^T |\phi_{11} - \phi_{21}|^2
  + |\phi_{12} - \phi_{22}|^2 \, d\mu(t) \Big)^{\frac12}\Big\|_{E}\\
  & \eqsim_{q} L_b\, (\| \phi_{11} -\phi_{21}\|_{\g(L^2((0,T),\mu),E)}
      + \| \phi_{12}  -\phi_{22}\|_{\g(L^2((0,T),\mu),E)})\\
      & \lesssim L_b \,\|\phi_{1} -\phi_{2}\|_{\g(L^2((0,T),\mu),X)}
\end{align*}
for all $\omega\in\O$.
The other estimate in (H4) can be  established  in a similar way.
\end{proof}

We next consider an $L^r(S)$--valued Brownian motion $w_1$, where $r\in [1, \infty]$.
In this case, we let $H_1$ be the
reproducing kernel Hilbert space of the Gaussian random variable $w_1(1,\cdot)$
and let $i_1$ be the embedding of $H_1$ into $L^r(S)$. Then  we have
$i_1 \in \g(H_1,L^r(S))$ and  $i_1 i_1^*\in \calL(L^{r'}(S),L^r(S))$
is the covariance operator of $w_1(1,\cdot)$ (cf.\ \cite{Bog}, \cite{BrzvNe}
and the references therein for details).
Let $W_{H_1}$ be a cylindrical Brownian motion such that $w_1 = i_1 W_{H_1}$.

\begin{example}\label{ex:Lrnoise}
Assume that $w_1$ is an $L^r(S)$-valued Brownian motion with $r>d$. Then
\ref{as:LipschitzBtype} is satisfied for all $\theta_B\in (\frac{d}{2r},\frac12)$
and $a=0$.
\end{example}

\begin{remark}
If $Q_1$ is of the form \eqref{eq:Q1} for an orthonormal system
$(e_n)_{n\ge1}$ in $L^2(S)$ and  numbers $\lambda_n\ge0$, then a
sufficient condition for $w_1$ to be $L^r$-valued is
\[\sum_{n\geq 1} \lambda_n \|e_n\|_{L^r(S)}^2<\infty,\]
or more generally
\begin{equation}\label{eq:assQ2}
\Big(\sum_{n\geq 1} \lambda_n |e_n|^2\Big)^{\frac12}\in L^r(S),
\end{equation}
cf.\ \cite{BrzLi}.
\end{remark}

\begin{proof}[Proof of Example~\ref{ex:Lrnoise}]
We use the same notation as in Example~\ref{ex:infty}, but $H_1$
will be the reproducing kernel Hilbert space for $w_1(1, \cdot)$ and
$(h_n)_{n\geq 1}$ is an orthonormal basis of $H_1$. By \cite[Proposition~2.1]{BP99} and \cite[Theorem 2.3]{BrzvN03} (also see \cite[Lemma~2.1]{NVW3})
\[\Big\|\Big(\sum_{n\geq 1} |i_1 h_n|^2 \Big)^{\frac12}\Big\|_{L^r(S)}
   \eqsim_r \|i_1\|_{\g(H_1,L^r(S))} <\infty.\]
We set $\tilde{\theta}_B = \frac{d}{2r}<\theta_B$ and
choose $\e>0$ with $\tilde{\theta}_B+\e<\theta_B$. Since
$E_{\frac12 - \frac12 (\tilde{\theta}_B+\e)}\times
E_{-\frac12(\tilde{\theta}_B+\e)} \hookrightarrow X_{-\theta_B}$
by \eqref{eq:inter/extra-damped2b}, we can estimate
\[
\begin{aligned}
\| B(\cdot,\omega,\phi_1)& -
B(\cdot,\omega,\phi_2)\|_{\g(L^2((0,T),\mu;H),X_{-\theta_B}))}
\\ & \lesssim_{\theta_B,r,d,q} \| B_2(\cdot,\omega,\phi_1) -
B_2(\cdot,\omega,\phi_2)\|_{\g(L^2((0,T),\mu;H),E_{-\frac12(\tilde{\theta}_B+\e)})}
\end{aligned}
\]
for each $\omega\in\O$.
Here $B_2$ is the second coordinate of $B$. The other one is zero.

Let $\frac1v=\frac1q+\frac1r$.
We claim that $L^v(S) \hookrightarrow E_{-\frac12(\tilde{\theta}_B+\e)}$.
Indeed, let $\B_v$ denote the realization of the negative Dirichlet
Laplacian in $L^v(S)$.
Taking into account Theorem~V.1.4.12 of \cite{Am}, we have to show that
\[\|x\|_{L^q(S)} \lesssim_{q,v,\tilde{\theta}_B,\e}
        \|\B_v^{\tilde{\theta}_B+\e} x\|_{L^v(S)} \]
for all $x\in D(\B_v^{\tilde{\theta}+\e})$. From
\cite[Theorem~4.3.1.2]{Tr1} and \eqref{eq:int-emb} we deduce
\begin{align*}
\|x\|_{B^{2\tilde{\theta}_B}_{v,1}(S)}
 &\eqsim_{\tilde{\theta}_B,v}\|x\|_{(L^v(S),W^{2,v}(S))_{\tilde{\theta}_B,1}}
 \eqsim_v \|x\|_{(L^v(S),D(B_v))_{\tilde{\theta}_B,1}} \\
 &\lesssim_{q,v,\tilde{\theta_B},\e}\|\B_v^{\tilde{\theta}_B+\e}x\|_{L^v(S)},
\end{align*}
so that the claim follows from Sobolev's embedding (cf.\
\cite[Theorem~4.6.1]{Tr1}). The claim, \eqref{eq:square}, H\"older's inequality
and the Lipschitz continuity of $b$ imply that
\begin{align*}
\| &B_2(\cdot,\omega,\phi_1) -B_2(\cdot,\omega,\phi_2)\|_{\g(L^2((0,T),\mu;H),E_{-\frac12({\tilde{\theta}_B}+\varepsilon)})}
\\ & \lesssim_{\theta_B,r,d,q} \| B_2(\cdot,\omega,\phi_1) - B_2(\cdot,\omega,\phi_2)\|_{\g(L^2((0,T),\mu;H),L^v(S))}
\\ & \eqsim_v \Big\|\Big( \int_0^T |b(t,\omega, \phi_{11},\phi_{12}) - b(t,\omega, \phi_{21}, \phi_{22})|^2 \, d\mu(t) \sum_{n\geq 1} |i_1 h_n|^2 \Big)^{\frac12}\Big\|_{L^v(S)}
\\ & \leq \Big\|\Big( \int_0^T \! |b(t,\omega, \phi_{11},\phi_{12}) - b(t,\omega, \phi_{21}, \phi_{22})|^2 \, d\mu(t) \Big)^{\frac12}\Big\|_{E}  \, \Big\|\Big(\sum_{n\geq 1} |i_1 h_n|^2 \Big)^{\frac12}\Big\|_{L^r(S)}
\\ & \lesssim_{w_1} L_b( \| \phi_{11} -\phi_{21}\|_{\g(L^2((0,T),\mu),E)} + \| \phi_{12}
-\phi_{22}\|_{\g(L^2((0,T),\mu),E)}).
\end{align*}
The other estimate in (H4) can be  established  in a similar way.
\end{proof}

If $d=1$ we can consider the space-time white noise
situation, where the covariance operator $Q_1:H_1\to H_1$
is the identity. This is possible
since in this case we can choose $q<d/(d-1)$ as large as needed.

\begin{example}\label{ex:d=1}
Let $Q_1= I$ on $H_1=L^2(S)$, $d=1$, and $q\in (2, \infty)$. Then
\ref{as:LipschitzBtype} is satisfied for all
$\theta_B>\frac{1}{4}+\frac{1}{2q}$ and $a=0$.
\end{example}

\begin{proof}
Let $q\in (2, \infty)$ and  $\frac12>\theta_B>\frac{1}{4}+\frac{1}{2q}$.
We take $\e>0$ be such that
$\theta_B-\e>\frac{1}{4}+\frac{1}{2q}$ and write $\theta_B -\e=
\theta_1+\theta_2$, where $\theta_1>\frac14$ and $\theta_2>\frac{1}{2q}$.
Since $L^q$ with $q\in (2, \infty)$ has type $2$,
Lemma~\ref{lem:type2Lipschitz} says that $B $ is $\gL$--Lipschitz
and of linear growth
if $B(t,\omega,\cdot):X\to \g(H_1,X_{-\theta_B})$ is Lipschitz
and of linear growth with a uniform constant.

We observe that $\mathcal{A}^{-\frac{\theta_1}{2}}\in \B(H_1,
W^{2\theta_1,2}(S))$ and that the injection $i:W^{2\theta_1,2}(S)\to L^q(S)$
belongs to $\g(W^{2\theta_1,2}(S),L^q(S))$ because of
 \cite[Corollary~2.2]{NVW3}. The right-ideal property \eqref{eq:ideal}
thus implies that
\[\|i\mathcal{A}^{-\frac{\theta_1}{2}}\|_{\g(H_1,L^q(S))}
\leq   \|i\|_{\g(W^{2\theta_1,2}(S),L^q(S))} \,
     \|\mathcal{A}^{-\frac{\theta_1}{2}}\|_{\calL(H_1,W^{2\theta_1,2}(S))}
<\infty.\]
For $x=(x_1,x_2)$ and $y=(y_1,y_2)$ in $X$, we deduce from
\eqref{eq:inter/extra-damped2b} and the right-ideal property that
\begin{align*}
\| B(t,x)&-B(t,y)\|_{\g(H_1,X_{-\theta_B})}   \lesssim_{\theta_B,q}
\| b(t,\omega,x_1, x_2)-b(t,\omega,y_1, y_2))\|_{\g(H_1,E_{-\frac12(\theta_B-\e)})}\\
& = \| i \mathcal{A}^{-\frac{\theta_1}{2}} \mathcal{A}^{-\frac{\theta_2}{2}}
( b(t,\omega, x_1, x_2)-b(t,\omega,y_1, y_2))\|_{\g(H_1,L^q(S))}\\
& \leq \| i \mathcal{A}^{-\frac{\theta_1}{2}}\|_{\g(H_1,L^q(S))}
   \|\mathcal{A}^{-\frac{\theta_2}{2}}
  (b(t,\omega, x_1, x_2)-b(t,\omega,y_1,y_2))\|_{\calL(H_1)}
\end{align*}
for all $\omega\in\O$ and $t\ge0$.
As in the claim in the proof of Example~\ref{ex:Lrnoise} one can use
Sobolev's embedding theorem to obtain
\begin{align*}
\|\mathcal{A}^{-\frac{\theta_2}{2}} (b(t,\omega, x_1, x_2)-&b(t,\cdot,y_1,y_2))\|_{\calL(H)}  \\
& \leq \|\mathcal{A}^{-\frac{\theta_2}{2}} (b(t,\cdot,
x_1, x_2)-b(t,\cdot,y_1, y_2))\|_{L^\infty(S)} \\
& \lesssim_{\theta_2,q} \|b(t,\cdot, x_1, x_2)-b(t,\cdot,y_1, y_2)\|_{E}\\
& \leq L_b \,(\|x_1-y_1\|_{E} +\|x_2-y_2\|_{E})   \lesssim L_b \|x-y\|_X.
\end{align*}
Thus we have shown the Lipschitz estimate in Lemma~\ref{lem:type2Lipschitz}.
The other estimate in this lemma  can be established in a similar way.
\end{proof}

\begin{remark}\label{rem:nabla}
It is clear from the proofs of Examples
\ref{ex:infty}, \ref{ex:Lrnoise} and \ref{ex:d=1}
that (H4) also holds if
$b$ also depends on $\nabla u$ and $\nabla^2 u$ in an appropriate
Lipschitz sense.  The same is true for $f$, $G$ and  $C$
in Theorem~\ref{thm:platen2}.
\end{remark}

\begin{remark}\label{rem:local}
In the above examples one could allow $f$ and $b$ to be only locally Lipschitz
in the third coordinate, i.e., the coordinate for $u(t,s)$. For this one needs
to define the maps $F,B,C,G$ on $X_a$ for a suitable $a>0$ such that the first
component of $X_a$ is embedded into $C(\overline{S})$. (See \cite[Theorems~8.1
and 10.2]{NVW3} for details.) This gives the condition $2+2a - \frac{d}{q}>0$.
However, we can only take $a>0$ such that $a+\theta_C<\frac12$ and
 $a+\theta_B<\frac12$. Since $\theta_C\in (\frac{d}{2 q'},\frac12)$ as
explained before Theorem~\ref{thm:platen2}, we obtain the first condition
$-1+\frac{d}{2}<\frac12$. This inequality holds for $d=1, 2$.

For Example~\ref{ex:infty} there are no conditions on $\theta_B$, so that
$d=1,2$ are both allowed. For Example~\ref{ex:Lrnoise} we also need
$\theta_B>\frac{d}{2r}$, and therefore $\frac{d}{2r} +
\frac{d}{2q}<\frac32$ must hold as well. This condition
holds for $d=1,2$ and all $r>d$ and $1<q<d/(d-1)$.
 For the Example~\ref{ex:d=1} we have $d=1$. There the condition reads
 $\theta_B >\frac14+\frac{1}{2q}$. Therefore, we obtain $\frac14 + \frac1q
<\frac12$. This holds if and only if $q>4$.
\end{remark}

\section{The damped wave equation\label{sec:stochwave}}

In this section we obtain existence, uniqueness and regularity results for a
structurally damped wave equation. Since the proofs follow the line
of arguments of the previous section, we omit the details.
The equation is given by
\begin{equation}\label{eq:wave}
\left\{
\begin{aligned}
\ddot{u}(t,s) - & \Delta u(t,s) - \rho (-\Delta)^{\frac12} \dot{u}(t,s)= f(t,s, u(t,s),\dot{u}(t,s)) \\
& + b(t,s,u(t,s),\dot{u}(t,s)) \, \frac{\partial w_1(t,s)}{\partial t} +
\Big[G(t,u(t,\cdot), \dot{u}(t,\cdot)) \\ & + C(t,u(t,\cdot), \dot{u}(t,\cdot))
\frac{\partial w_2(t)}{\partial t} \Big] \delta(s-s_0), \ \ t\in [0,T], s\in S,
\\ & u(0,s) = u_0(s), \ \dot{u}(0,s) = u_1(s), \ \ s\in S,
\\ & u(t,s) = 0,\ \ t\in[0,T], \  s\in \partial S,
\end{aligned}\right.
\end{equation}
where $S\subset \R^n$ has a $C^2$ boundary $\partial S$
and $(-\Delta)^{\frac12}$ denotes the square root of the negative
Dirichlet Laplacian. We reformulate this equation as \eqref{eq:SEtype}
in the same way as in Section~\ref{sec:plate}.

Let $q\in (1, \infty)$ and  $E = L^q(S)$. On $E$ we define
$(\mathcal{A},D(\mathcal{A}))$ by
\[\mathcal{A} x = -\Delta x, \ \ D(\mathcal{A})  = W^{2,q}(S)\cap W^{1,q}_0(S).\]
Let $X = E_{\frac12}\times E$ and define $(A,D(A))$ by
\[A = \left(%
\begin{array}{cc}
  0 & I \\
  \mathcal{A} & -\rho \mathcal{A}^{\frac12} \\
\end{array}%
\right), \ \ D(A) = D(\mathcal{A})\times D(\mathcal{A}^{\frac12}).
\]
It follows from Proposition~\ref{prop:damped} that $A$ generates an
analytic semigroup $(S(t))_{t\geq 0}$.

We further assume that $\rho>0$ and $s_0\in S$ are fixed and that
$\delta$ is the usual point evaluation. Moreover,
$f$, $b$, $C$, $G$, $w_1$, $w_2$ shall satisfy the assumptions (A0)--(A3)
 in Section~\ref{sec:plate} for the above space $X$ and the maps
 $F$, $B$ and $\Lambda$ are defined  as in Section~\ref{sec:plate} for
the above space $X$. Finally, it assumed that $B$ fulfills
hypothesis (H4). Noting that $\mathcal{A}$ is now of second order,
one can see in the same way as in Section~\ref{sec:plate} that
$\Lambda$ is well-defined for all $1<q<\frac{2d}{2d-1}$. A mild and
weak solution are defined in a similar way as in Section
\ref{sec:plate}.
Finally, for $1<q<2$ we have
\[E_{\frac12 +\frac12 \delta}\times E_{\frac12 \delta} = (W^{1+\delta,q}(S)\cap W^{1,q}_0(S)) \times W^{\delta,q}(S).\]

\begin{theorem}\label{thm:wave}
Let $1<q<\frac{2d}{2d-1}$. Assume that $u_0:\O\to W^{1,q}_0(S)$ and
$u_1:\O\to L^q(S)$ are $\F_0$-measurable. Let $f,G,b,C,w_1$ and
$w_2$ be as above. The following assertions hold:

\begin{enumerate}
\item For all $\alpha\in (0,\frac12)$ and $p>2$ such that
$a+\max\{\theta_B,\theta_C\}<\alpha-\frac1p$, there exists a unique
mild solution $U$ of \eqref{eq:wave} in
\[V^{0}_{\a,p}([0,T]\times\O;W^{1,q}_0(S) \times L^q(S)).\]
There is a unique weak solution $u\in W^{1,2}(0,T;L^q(S))$ of
\eqref{eq:plate} such that $(u,\dot{u})$ belongs to the space in
\eqref{eq:solspace}. Moreover, $U = (u, \dot{u})$.

\item There exists a version of $u$ with paths that satisfy $u\in
C([0,T];W^{1,q}_0(S))$ and $\dot{u} \in C([0,T];L^q(S))$.

\item Let $\eta\in (0,\frac12]$. If $u_0:\O\to H^{1+\eta,q}(S)\cap W^{1,q}_0(S)$
and $u_1:\O\to H^{\eta,q}(S)$, then there exists a version of $u$ with paths that
satisfy $u\in C^{\lambda}([0,T];H^{1+\delta,q}(S)\cap W^{1,q}_0(S))$ and
$\dot{u} \in C^{\lambda}([0,T];H^{\delta,q}(S))$ for all $\delta,\lambda\geq 0$
such that $\delta+\lambda<\min\{\eta, \frac12-\theta_B, \frac12-\theta_C\}$.
\end{enumerate}
\end{theorem}

This theorem can be proved in the same way as
Theorem~\ref{thm:platen2}. Let us give some examples for $w_1$.
Example~\ref{ex:infty} works in exactly the same way for the wave
equation. Example~\ref{ex:Lrnoise} has the following version for the
wave equation.

\begin{example}
Assume that $w_1$ is an $L^r(S)$-valued Brownian motion with $r>2d$.
Then \ref{as:LipschitzBtype} is satisfied for all $\theta_B\in
(\frac{d}{r}, 1)$ and $a=0$.
\end{example}

\noindent
This assertion can be shown as in Example~\ref{ex:Lrnoise},
we thus leave the details to reader.

\def\cprime{$'$} \def\polhk#1{\setbox0=\hbox{#1}{\ooalign{\hidewidth
  \lower1.5ex\hbox{`}\hidewidth\crcr\unhbox0}}} \def\cprime{$'$}
\providecommand{\bysame}{\leavevmode\hbox to3em{\hrulefill}\thinspace}

\end{document}